\documentstyle[12pt]{article}
\begin{document}
\author{S.V. Ludkovsky.}
\title{Smoothness of functions global and along curves
over ultra-metric fields.}
\date{23.08.2006}
\maketitle

\begin{abstract}
The article is devoted to the investigation of smoothness of
functions $f(x_1,...,x_m)$ of variables $x_1,...,x_m$ in infinite
fields with non trivial non archimedean valuations, where $m\ge 2$.
Theorems about classes of smoothness $C^n$ or $C^n_b$ of functions
with continuous or bounded uniformly continuous on bounded domains
partial difference quotients up to the order $n$ are investigated.
It is proved that from $f\circ u\in C^n({\bf K},{\bf K}^l)$ or
$f\circ u\in C^n_b({\bf K},{\bf K}^l)$ for each $C^{\infty }$ or
$C^{\infty }_b$ curve $u: {\bf K}\to {\bf K}^m$ it follows, that
$f\in C^n({\bf K}^m,{\bf K}^l)$ or $f\in C^n_b({\bf K}^m,{\bf
K}^l)$. Moreover, classes of smoothness $C^{n,r}$ and $C^{n,r}_b$
and more general in the sense of Lipschitz for partial difference
quotients are considered and theorems for them are proved.
\end{abstract}

\section{Introduction}
Fields with non archimedean valuations such as the field of $p$-adic
numbers were first introduced by K. Hensel \cite{hensel}. Then it
was proved by A. Ostrowski \cite{ostrow} that on the field of
rational numbers each multiplicative norm is either the usual norm
as in $\bf R$ or is equivalent to a non archimedean norm
$|x|=p^{-k}$, where $x=np^k/m\in \bf Q$, $n, m, k\in \bf Z$, $p\ge
2$ is a prime number, $n$ and $m$ and $p$ are mutually pairwise
prime numbers. It is well known, that each locally compact infinite
field with a non trivial non archimedean valuation is either a
finite algebraic extension of the field of $p$-adic numbers or is
isomorphic to the field ${\bf F}_{p^k}(\theta )$ of power series of
the variable $\theta $ with expansion coefficients in the finite
field ${\bf F}_{p^k}$ of $p^k$ elements, where $p\ge 2$ is a prime
number, $k\in \bf N$ is a natural number \cite{roo,weil}. Non
locally compact fields are also wide spread \cite{diar,roo,sch1}.
\par Last years non archimedean analysis \cite{roo,sch1,sch2}
and mathematical physics \cite{yojan,khrhs,khrqp,vla3} are being
fastly delevoped. But many questions and problems remain open. In
the classical case it is known the Boman's theorem relating
smoothness of a function of several real variables and of its
compositions with smooth curves \cite{boman,krimich}. But this
problem was not studied completely in the non archimedean case
besides particular cases and instead of curves for compositions with
functions of more than one variable \cite{beglneeb}, that is the
significant simplification of the problem.
\par In the non archimedean analysis classes of smoothness are defined
in another fashion as in the classical case over $\bf R$, since
locally constant functions on fields $\bf K$ with non archimedean
valuations are infinite differentiable and there exist non trivial
non locally constant functions infinite differentiable with
identically zero derivatives \cite{sch1,sch2}. This is caused by the
stronger ultra-metric inequality $|x+y|\le \max (|x|,|y|)$ in
comparison with the usual triangle inequality, where $|x|$ is a
multiplicative norm in $\bf K$ \cite{roo}. In papers
\cite{beglneeb,luseamb,lutmf99,luanmbp} there were considered
classes of smoothness $C^n$ for functions of several variables in
non archimedean fields or in topological vector spaces over such
fields.
\par This paper is devoted to
the investigation of smoothness of functions $f(x_1,...,x_m)$ of
variables $x_1,...,x_m$ in infinite fields with non trivial non
archimedean valuations, where $m\ge 2$. In the paper fields locally
compact and as well as non locally compact are considered. Theorems
about classes of smoothness $C^n$ or $C^n_b$ of functions with
continuous or uniformly continuous on bounded domains partial
difference quotients up to the order $n$ are investigated. It is
proved that from $f\circ u\in C^n({\bf K},{\bf K})$ or $f\circ u\in
C^n_b({\bf K},{\bf K})$ for each $C^{\infty }$ or $C^{\infty }_b$
curve $u$ it follows, that $f\in C^n({\bf K}^m,{\bf K})$ or $f\in
C^n_b({\bf K}^m,{\bf K})$. Moreover, classes of smoothness $C^{n,r}$
and $C^{n,r}_b$ and more general in the sense of Lipschitz for
partial difference quotients are considered and theorems for them
are proved. \par Many specific features of the non archimedean case
in comparison with the classical one are found. In the non
archimedean case analogs of classical theorems over $\bf R$ such as
3 and 10 \cite{boman} are not true due to the ultrametric inequality
for the non archimedean norm, and since if a function $f$ is
homogeneous, then ${\bar {\Phi }}^k$ need not be homogeneous for
$k\ge 1$. Theorem 2 from \cite{boman} in the non archimedean case is
true in the stronger form due to the ultrametric inequality (see
Theorem 39 below). The notion of quasi analyticity used in the
classical case in \cite{boman} has not sense in the non archimedean
case because of the necessity to operate with the partial difference
quotients ${\bar {\Phi }}^kf$ instead of derivatives $D^kf$. It
leads naturally to the local analyticity in the non archimedean
case. In the latter case the exponential function has finite radius
of convergence on $\bf K$ with $char ({\bf K})=0$. Moreover, in the
proof of Theorem 42 it was used specific feature of the non
archimedean analysis of analytic functions for which an analog of
the Louiville theorem is not true (see also \cite{sch1}).
\par Several lemmas of the paper serve for subsequent proofs
of theorems. It is proved in theorems 38-42 below, that for
corresponding smoothness, for example, $C^n_{\phi }({\bf K}^m,{\bf
K})$ of a function $f$ it is sufficient that $f\circ u\in C^n_{\phi
} ({\bf K},{\bf K})$ for each curve $u\in C^{\infty }({\bf K},{\bf
K}^m)$, but a local analyticity of $u$ instead of $C^{\infty }$ is
insufficient.

\section{Smoothness of functions}
\par {\bf 1. Definitions.} Let $\bf K$ be an infinite field with a
non trivial non archimedean valuation, let also $X$ and $Y$ be
topological vector spaces over $\bf K$ and $U$ be an open subset in
$X$. For a function $f: U\to Y$ consider the associated function
\par $f^{[1]}(x,v,t) := [f(x+tv) - f(x)]/t$ \\
on a set $U^{[1]}$ at first for $t\ne 0$ such that $U^{[1]} := \{
(x,v,t)\in X^2\times {\bf K}, x\in U, x+tv\in U \} $. If $f$ is
continuous on $U$ and $f^{[1]}$ has a continuous extension on
$U^{[1]}$, then we say, that $f$ is continuously differentiable or
belongs to the class $C^1$. The $\bf K$-linear space of all such
continuously differentiable functions $f$ on $U$ is denoted
$C^{[1]}(U,Y)$. By induction we define functions $f^{[n+1]}:=
(f^{[n]})^{[1]}$ and spaces $C^{[n+1]}(U,Y)$ for $n=1,2,3,...$,
where $f^{[0]}:=f$, $f^{[n+1]}\in C^{[n+1]}(U,Y)$ has as the domain
$U^{[n+1]} := (U^{[n]})^{[1]}$.
\par The differential $df(x): X\to Y$ is defined as
$df(x)v := f^{[1]}(x,v,0)$.
\par Define also partial difference quotient operators $\Phi ^n$
by variables corresponding to $x$ only such that
\par $\Phi ^1f(x;v;t) = f^{[1]}(x,v,t)$ \\
at first for $t\ne 0$ and if $\Phi ^1f$ is continuous for $t\ne 0$
and has a continuous extension on $U^{[1]}=:U^{(1)}$, then we denote
it by ${\bar {\Phi }}^1f(x;v;t)$. Define by induction \par $\Phi
^{n+1} f(x;v_1,...,v_{n+1};t_1,...,t_{n+1}):= \Phi ^1(\Phi
^nf(x;v_1,...,v_n;t_1,...,t_n))(x;v_{n+1};t_{n+1})$ \\
at first for $t_1\ne 0,...,t_{n+1}\ne 0$ on $U^{(n+1)}:= \{
(x;v_1,...,v_{n+1};t_1,...,t_{n+1}): x\in U; v_1,...,v_{n+1}\in X;
t_1,...,t_{n+1}\in {\bf K}; x+v_1t_1\in
U,...,x+v_1t_1+...+v_{n+1}t_{n+1}\in U \} $. If $f$ is continuous on
$U$ and partial difference quotients $\Phi ^1f$,...,$\Phi ^{n+1}f$
has continuous extensions denoted by ${\bar {\Phi }}^1f$,..., ${\bar
{\Phi }}^{n+1}f$ on $U^{(1)}$,...,$U^{(n+1)}$ respectively, then we
say that $f$ is of class of smoothness $C^{n+1}$. The $\bf K$ linear
space of all $C^{n+1}$ functions on $U$ is denoted by
$C^{n+1}(U,Y)$, where $\Phi ^0f := f$, $C^0(U,Y)$ is the space of
all continuous functions $f: U\to Y$. Then the differential is given
by the equation $d^nf(x).(v_1,...,v_n) := n! {\bar {\Phi
}}^nf(x;v_1,...,v_n;0,...,0)$, where $n\ge 1$, also denote
$D^nf=d^nf$. Shortly we shall write the argument of $f^{[n]}$ as
$x^{[n]}\in U^{[n]}$ and of ${\bar {\Phi }}^nf$ as $x^{(n)}\in
U^{(n)}$, where $x^{[0]}=x^{(0)}=x$, $x^{[1]}=x^{(1)}=(x,v,t)$,
$v^{[0]}=v^{(0)}=v$, $t_1=t$, $x^{[k]}=(x^{[k-1]},v^{[k-1]},t_k)$
for each $k\ge 1$, $x^{(k)} := (x;v_1,...,v_k;t_1,...,t_k)$.
\par Subspaces of $C^n$ or $C^{[n]}$ of all bounded uniformly continuous
functions together with ${\bar {\Phi }}^kf$ or $\Upsilon ^kf$ on
bounded open subsets of $U$ and $U^{(k)}$ or $U^{[k]}$ for
$k=1,...,n$ denote by $C^n_b(U,Y)$ or $C^{[n]}_b(U,Y)$ respectively.
\par Consider partial difference quotients of products and compositions
of functions and relations between partial difference quotients and
differentiability of both types. Denote by $L(X,Y)$ the space of all
continuous $\bf K$-linear mappings $A: X\to Y$. By $L_n(X^{\otimes
n},Y)$ denote the space of all continuous $\bf K$ $n$-linear
mappings $A: X^{\otimes n}\to Y$, particularly,
$L(X,Y)=L_1(X^{\otimes 1},Y)$. If $X$ and $Y$ are normed spaces,
then $L_n(X^{\otimes n},Y)$ is supplied with the operator norm: $ \|
A \| := \sup_{h_1\ne 0,...,h_n\ne 0; h_1,...,h_n\in X} \|
A.(h_1,...,h_n) \|_Y/ (\| h_1 \| _X... \| h_n \| _X)$.
\par {\bf 2. Lemma.} {\it The spaces $C^{[1]}(U,Y)$ and $C^1(U,Y)$
are linearly topologically isomorphic. If $f\in C^n(U,Y)$, then
${\bar {\Phi }}^nf(x;*;0,...,0): X^{\otimes n}\to Y$ is a $\bf K$
$n$-linear $C^0(U,L_n(X^{\otimes n},Y))$ symmetric map.}
\par {\bf Proof.} From Definition 1 it follows, that
$f^{[1]}(x,v,t)= {\bar {\Phi }}^1f(x;v;t)$ on $U^{[1]}=U^{(1)}$, so
both $\bf K$-linear spaces are linearly topologically isomorphic. On
the other hand, it was proved in Proposition 2.2 and Lemma 4.8
\cite{beglneeb}, that ${\bar {\Phi }}^nf(x;*,0,...,0)$ is the $\bf
K$ $n$-linear symmetric mapping for each $x\in U$ and it belongs to
$C^0(U,L_n(X^{\otimes n},Y))$, since ${\bar {\Phi
}}^nf(x;v_1,...,v_n;t_1,...,t_n)$ is continuous on $U^{(n)}$ and for
each $x\in U$ and $v_1,...,v_n\in X$ there exist neighborhoods $V_i$
of $v_i$ in $X$ and $W$ of zero in $\bf K$ such that
$x+WV_1+...+WV_n\subset U$.
\par {\bf 3. Lemma.} {\it Operators $\Upsilon ^n(f) := f^{[n]}$
from $C^{[n]}(U,Y)$ into $C^0(U^{[n]},Y)$ and ${\bar {\Phi }}^n:
C^n(U,Y)\to C^0(U^{(n)},Y)$ are $\bf K$-linear and continuous.}
\par {\bf Proof.} Since $[(af+bg)(x+vt)-(af+bg)(x)]/t=
a(f(x+vt)-f(x))/t + b(g(x+vt)-g(x))/t$ for each $f, g\in C^1(U,Y)$
and each $a, b\in \bf K$, then applying this formula by induction
and using definitions of operators $\Upsilon ^n$ and ${\bar {\Phi
}}^n$ we get their $\bf K$-linearity. Indeed, \par $\Upsilon
^n(af+bg)(x^{[n]})=\Upsilon ^1(\Upsilon
^{n-1}(af+bg)(x^{[n-1]}))(x^{[n]})=\Upsilon
^1(af^{[n-1]}+bg^{[n-1]})(x^{[n]})=af^{[n]}(x^{[n]})+bg^{[n]}(x^{[n]})$
and
\par ${\bar {\Phi }}^n(af+bg)(x^{(n)})={\bar {\Phi }}^1(
{\bar {\Phi }}^{n-1}(af+bg)(x^{(n-1)}))(x^{(n)})={\bar {\Phi
}}^1(af^{(n-1)}+bg^{(n-1)})(x^{(n)})=af^{(n)}(x^{(n)})+bg^{(n)}(x^{(n)})$.
\\ The continuity of $\Upsilon ^n$ and ${\bar {\Phi }}^n$ follows
from definitions of spaces $C^{[n]}(U,Y)$ and $C^n(U,Y)$
respectively.
\par {\bf 4. Lemma.} {\it  Let either
$f, g\in C^{[n]}(U,Y)$, where $U$ is an open subset in $X$, $Y$ is
an algebra over $\bf K$, or $f\in C^{[n]}(U,{\bf K})$ and $g\in
C^{[n]}(U,Y)$, where $Y$ is a topological vector space over $\bf K$,
then
\par $(1)$ $(fg)^{[n]}(x^{[n]}) = ({\Upsilon }\otimes {\hat P} +
{\hat {\pi }}\otimes {\Upsilon })^n
.(f\otimes g)(x^{[n]})$ \\
and $(fg)^{[n]}\in C^0(U^{[n]},Y)$, where $({\hat {\pi
}}^kg)(x^{[k]}) := g\circ \pi ^0_1\circ \pi ^1_2\circ ...\circ \pi
^{k-1}_k(x^{[k]})$, ${\hat P}^ng := P_nP_{n-1}...P_1g$, $\pi
^{k-1}_k(x^{[k]}) := x^{[k-1]}$, $(A\otimes B).(f\otimes g):=
(Af)(Bg)$ for $A, B\in L(C^n(U,Y),C^m(U,Y))$, $m\le n$, $(A_1\otimes
B_1)...(A_k\otimes B_k).(f\otimes g):=(A_1...A_k\otimes
B_1...B_k).(f\otimes g):= (A_1...A_kf)(B_1...B_kg)$ for
corresponding operators, ${\Upsilon
}^nf := f^{[n]}$, $(P_kg)(x^{[k]}) := g(x^{[k-1]}+v^{[k-1]}t_k)$, \\
${\hat P}^k{\hat {\pi }}^{a_1}\Upsilon ^{b_1}...{\hat {\pi
}}^{a_l}\Upsilon ^{b_l}g = P_{k+s} ...P_{s+1}{\hat {\pi
}}^{a_1}\Upsilon ^{b_1}...{\hat {\pi }}^{a_l}\Upsilon ^{b_l}g$ with
$s=b_1+...+b_l-a_1-...-a_l\ge 0$, $a_1,...,a_l, b_1,...,b_l \in \{
0, 1, 2, 3,... \} $.}
\par {\bf Proof.} Let at first $n=1$, then
\par $(2)$ $(fg)^{[1]}(x^{[1]})=[(fg)(x+vt)-(fg)(x)]/t =
[(f(x+vt)-f(x))g(x+vt) + f(x)(g(x+vt)-g(x))]/t= ({\Upsilon
}^1f)(x^{[1]})(P_1g)(x^{[1]}) + ({\hat {\pi
}}^0_1f)(x^{[1]}){\Upsilon }^1g(x^{[1]})$, \\
since ${\hat {\pi }}^0_1(x^{[1]})=x$ and $P_1$ is the composition of
the projection ${\hat {\pi }}^0_1$ and the shift operator on $vt$.
Let now $n=2$, then applying Formula $(2)$ we get: \par $(3)$
$(fg)^{[2]}(x^{[2]})= ((fg)^{[1]}(x^{[1]}))^{[1]}(x^{[2]}) =
(\Upsilon ^1(f^{[1]}(x^{[1]})(x^{[2]}))g(x+(v^{[0]}+v^{[1]}_2t_2)
(t_1+v^{[1]}_3t_2) + v^{[1]}_1t_2) + f^{[1]}(x^{[1]})
g^{[1]}(x+v^{[0]}t_1,v^{[1]}_1+v^{[1]}_2 (t_1+v^{[1]}_3t_2),t_2) +
f^{[1]}(x,v^{[1]}_1,t_2) g^{[1]}(x^{[1]}+v^{[1]}_1t_2) + f(x)
g^{[2]}(x^{[2]})$, \\
where $v^{[k]} = (v^{[k]}_1, v^{[k]}_2, v^{[k]}_3)$ for each $k\ge
1$ and $v^{[0]}=v^{[0]}_1$ such that $x^{[k]}+v^{[k]}t_{k+1} =
(x^{[k]}+v^{[k]}_1t_{k+1}, v^{[k-1]}+ v^{[k]}_2t_{k+1},
t_k+v^{[k]}_3t_{k+1})$ for each $1\le k\in \bf Z$. For $n=3$ we get
\par $(4)$ $(fg)^{[3]}(x^{[3]})=[(\Upsilon ^3f)({\hat P}^3g) +
({\hat {\pi }}^1\Upsilon ^2f)(\Upsilon ^1{\hat P}^2g)+ (\Upsilon ^1
({\hat {\pi }}^1\Upsilon ^1f))({\hat P}^1{\Upsilon ^1}{\hat P}^1g)
+({\hat {\pi }}^2\Upsilon ^1f)(\Upsilon ^2{\hat P}^1g)+ (\Upsilon
^2{\hat {\pi }}^1f)({\hat P}^2\Upsilon ^1g)+ ({\hat {\pi
}}^1{\Upsilon }^1{\hat {\pi }}^1f) (\Upsilon ^1{\hat P}^1\Upsilon
^1g)+ (\Upsilon ^1({\hat {\pi }}^2f))({\hat
P}^1\Upsilon ^2g)+ ({\hat {\pi }}^3f) (\Upsilon ^3g)](x^{[3]})$, \\
since by our definition ${\hat P}^k{\hat {\pi }}^{a_1}\Upsilon
^{b_1}...{\hat {\pi }}^{a_l}\Upsilon ^{b_l}g = P_{k+s}
...P_{s+1}{\hat {\pi }}^{a_1}\Upsilon ^{b_1}...{\hat {\pi
}}^{a_l}\Upsilon ^{b_l}g$ with $s=b_1+...+b_l-a_1-...-a_l\ge 0$,
$a_1,...,a_l, b_1,...,b_l \in \{ 0, 1, 2, 3,... \} $. \par
Therefore, Formula $(1)$ for $n=1$ and $n=2$ and $n=3$ is
demonstrated by Formulas $(2-4)$. If $f, g \in C^0(U^{[k]},Y)$, $a,
b\in \bf K$, then $(P_k(af+bg))(x^{[k]}) :=
(af+bg)(x^{[k-1]}+v^{[k-1]}t_k)=$ $af(x^{[k-1]}+v^{[k-1]}t_k)+
bg(x^{[k-1]}+v^{[k-1]}t_k)$, moreover, ${\hat {\pi
}}^k(af+bg)(x^{[k]}) = (af+bg)\circ \pi ^0_1\circ \pi ^1_2\circ
...\circ \pi ^{k-1}_k(x^{[k]})=(af+bg)(x)= a f(x)+ b g(x) = a{\hat
{\pi }}^kf(x^{[k]}) + b {\hat {\pi }}^kg(x^{[k]})$ for each
$x^{[k]}\in U^{[k]}$, hence ${\hat {\pi }}^k$ and $P_k$ and ${\hat
P}^k$ are $\bf K$-linear operators for each $k\in \bf N$. Suppose
that Formula $(1)$ is proved for $n=1,...,m$, then for $n=m+1$ it
follows by application of Formula $(2)$ to both sides of Formula
$(1)$ for $n=m$:
\par $(fg)^{m+1}(x^{[m+1]})=((fg)^{[m]}(x^{[m]}))^{[1]}(x^{[m+1]})=
(({\Upsilon }\otimes {\hat P} + {\hat {\pi }}\otimes {\Upsilon })^m
.(f\otimes g)(x^{[m]}))^{[1]}(x^{[m+1]})= ({\Upsilon }\otimes {\hat
P} + {\hat {\pi }}\otimes {\Upsilon })^{m+1} .(f\otimes
g)(x^{[m+1]})$, \\
since $x^{[m+1]}=(x^{[m]})^{[1]}$ and more generally
$x^{[m+k]}=(x^{[m]})^{[k]}$ for each nonnegative integers $m$ and
$k$ such that $\pi ^{k-1}_k(x^{[m+k]})=x^{[m+k-1]}$ for $k\ge 1$;
$\Upsilon ^k$, ${\hat P}^k$ and ${\hat {\pi }}$ are $\bf K$-linear
operators on corresponding spaces of functions (see above and Lemma
3) and \par $({\Upsilon }\otimes {\hat P} + {\hat {\pi }}\otimes
{\Upsilon })^{m+1} .(f\otimes g)(x^{[m+1]})=$
\\  $\sum_{a_1+...+a_{m+1}+b_1+...+b_{m+1}=m+1} ({\Upsilon }^{a_1}\otimes
{\hat P}^{a_1})$\\ $({\hat {\pi }}^{b_1}\otimes {\Upsilon
}^{b_1})... ({\Upsilon }^{a_{m+1}}\otimes {\hat P}^{a_{m+1}}) ({\hat
{\pi }}^{b_{m+1}}\otimes {\Upsilon }^{b_{m+1}}).(f\otimes
g)(x^{[m+1]}) $, \\
where $a_j$ and $b_j$ are nonnegative integers for each
$j=1,...,m+1$, $(A_1\otimes B_1)...(A_k\otimes B_k).(f\otimes g):=
(A_1...A_k\otimes B_1...B_k).(f\otimes
g):=(A_1...A_kf)(B_1...B_kg).$
\par {\bf 5. Note.} Consider the projection
\par $(1)$ $\psi _n: X^{m(n)}\times {\bf K}^{s(n)}\to
X^{l(n)}\times {\bf K}^n$, \\
where $m(n)=2m(n-1)$, $s(n)=2s(n-1)+1$, $l(n)=n+1$ for each $n\in
\bf N$ such that $m(0)=1$, $s(0)=0$, $m(n)=2^n$,
$s(n)=1+2+2^2+...+2^{n-1}=2^n-1$. Then $m(n)$, $s(n)$, $l(n)$ and
$n$ correspond to number of variables in $X$, $\bf K$ for $\Upsilon
^n$, in $X$ and $\bf K$ for ${\bar {\Phi }}^n$ respectively.
Therefore, $\psi (x^{[n]})=x^{(n)}$ and $\psi _n(U^{[n]})=U^{(n)}$
for each $n\in \bf N$ for suitable ordering of variables. Thus
${\bar {\Phi }}^nf(x^{(n)})= {\hat {\psi }}_n{\Upsilon
}^nf(x^{[n]})=f^{[n]}(x^{[n]})|_{W^{(n)}}$, where ${\hat {\psi
}}_ng(y) := g(\psi _n(y))$ for a function $g$ on a subset $V$ in
$X^{l(n)}\times {\bf K}^n$ for each $y\in \psi _n^{-1}(V)\subset
X^{m(n)}\times {\bf K}^{s(n)}$, $W^{(n)}=U^{(n)} \times 0$, $0\in
X^{m(n)-l(n)}\times {\bf K}^{s(n)-n}$ for the corresponding ordering
of variables.
\par {\bf 6. Corollary.} {\it Let either $f, g\in C^n(U,Y)$, where $U$
is an open subset in $X$, $Y$ is an algebra over $\bf K$, or $f\in
C^n(U,{\bf K})$ and $g\in C^n(U,Y)$, where $Y$ is a topological
vector space over $\bf K$, then
\par $(1)$ ${\bar {\Phi }}^n(fg)(x^{(n)}) =
({\bar {\Phi }}\otimes {\hat P} + {\hat {\pi }}\otimes {\bar {\Phi }
})^n .(f\otimes g)(x^{(n)})$ \\
and ${\bar {\Phi }}^n(fg)\in C^0(U^{(n)},Y)$. In more details:
\par $(2)$ ${\bar {\Phi }}^n(fg)(x^{(n)}) =
\sum_{0\le a, 0\le b, a+b=n}\sum_{j_1<...<j_a; s_1<...<s_b; \{
j_1,...,j_a \} \cup \{ s_1,...,s_b \} = \{ 1,...,n \} }$ \\
${\bar {\Phi }}^af(x;v_{j_1},...,v_{j_a};t_{j_1},...,t_{j_a}) {\bar
{\Phi }}^bg(x+v_{j_1}t_{j_1}+...+v_{j_a}t_{j_a};v_{s_1},...,v_{s_b};
t_{s_1},...,t_{s_b})$.}
\par {\bf Proof.} The operator ${\hat {\psi }}_n$ is $\bf K$-linear,
since ${\hat {\psi }}_n(af+bg)(y)=(af+bg)(\psi _n(y))=af(\psi
_n(y))+ bg(\psi _n(y))$ for each $a, b\in \bf K$ and functions $f,
g$ on a subset $V$ in $X^{l(n)}\times {\bf K}^n$ and each $y\in \psi
_n^{-1}(V)\subset X^{m(n)}\times {\bf K}^{s(n)}$. Mention that the
restrictions of ${\hat {\pi }}^{k-1}_k$ and $P_k$ on $W^{(k)}$ gives
$\pi ^{k-1}_k(x^{(k)}) := x^{(k-1)}$ and $(P_kg)(x^{(k)}) :=
g(x^{(k-1)}+v_kt_k)$ in the notation of \S 1. The application of the
operator ${\hat {\psi }}_n$ to both sides of Equation $4(1)$ gives
Equation $(1)$ of this corollary, since ${\hat {\psi }}_n\Upsilon ^n
= {\bar {\Phi }}^n$ for each nonnegative integer $n$, where
$\Upsilon ^0=I$ and ${\bar {\Phi }}^0=I$ and ${\hat {\psi }}_0=I$
are the unit operators.
\par {\bf 7. Lemma.} {\it Let $f_1,...,f_k\in C^{[n]}(U,Y)$, where $U$
is an open subset in $X$, either $Y$ is an algebra over $\bf K$, or
$f_1,...,f_{k-1}\in C^{[n]}(U,{\bf K})$ and $f_k\in C^{[n]}(U,Y)$,
where $Y$ is a topological vector space over $\bf K$, then
\par $(1)$ $(f_1...f_k)^{[n]}(x^{[n]}) = [\sum_{\alpha =0}^{k-1}
{\hat {\pi }}^{\otimes \alpha }\otimes \Upsilon \otimes {\hat
P}^{\otimes (k-\alpha -1)}]^n.(f_1\otimes ...\otimes f_k)(x^{[n]})$ \\
and $(f_1...f_k)^{[n]}\in C^0(U^{[n]},Y)$, where ${\hat {\pi
}}^{\otimes \alpha }\otimes \Upsilon \otimes {\hat P}^{\otimes
(k-\alpha -1)}.(f_1\otimes ... \otimes f_k) := ({\hat {\pi
}}(f_1...f_{\alpha }))(\Upsilon f_{\alpha +1})({\hat P}(f_{\alpha
+2}...f_k))$, where ${\hat {\pi }}^0:=I$, ${\hat P}^0=I$ is the unit
operator, ${\hat {\pi }}f_0:=1$, ${\hat P}f_{k+1}:=1$ (see Lemma
4).}
\par {\bf Proof.} Consider at first $n=1$ and apply Formula $4(1)$
by induction to appearing products of functions, then
\par $(2)$ $\Upsilon ^1(f_1...f_k)(x^{[1]})=[(\Upsilon
^1(f_1...f_{k-1}))(P_1f_k) + ({\hat {\pi
}}^1(f_1...f_{k-1}))(\Upsilon ^1f_k)](x^{[1]})= [(\Upsilon
^1(f_1...f_{k-2}))(P_1f_{k-1})(P_1f_k) + ({\hat {\pi
}}^1(f_1...f_{k-2}))(\Upsilon ^1f_{k-1})(P_1f_k)$\\ $+ ({\hat {\pi
}}^1(f_1...f_{k-1}))(\Upsilon ^1f_k)](x^{[1]}) =...$\\
$=(\sum_{\alpha =0}^{k-1}({\hat {\pi }}^1)^{\otimes \alpha }\otimes
\Upsilon ^1\otimes P_1^{\otimes
(k-\alpha -1)}).(f_1\otimes ...\otimes f_k)$, \\
where $A^{\otimes \alpha }\otimes B\otimes C^{\otimes (k-\alpha
-1)}.(f_1\otimes ... \otimes f_k) := (A(f_1...f_{\alpha
}))(Bf_{\alpha +1})(C(f_{\alpha +2}...f_k))$ for operators $A, B$
and $C$ and each nonnegative integer $\alpha $, where $A^0:=I$,
$C^0=I$ is the unit operator, $Af_0:=1$, $Cf_{k+1}:=1$, in
particular, $A={\hat {\pi }}^1$, $B=\Upsilon ^1$, $C=P_1$. Thus,
acting by induction on both sides by $\Upsilon ^1$ from Formula
$(2)$ we get Formula $(1)$ of this lemma, since the product of $n$
terms $\Upsilon ^1...\Upsilon ^1$  is equal to $\Upsilon ^n$.
\par {\bf 8. Corollary.} {\it Let $f_1,...,f_k\in C^n(U,Y)$, where $U$
is an open subset in $X$, either $Y$ is an algebra over $\bf K$, or
$f_1,...,f_{k-1}\in C^n(U,{\bf K})$ and $f_k\in C^n(U,Y)$, where $Y$
is a topological vector space over $\bf K$, then
\par $(1)$ ${\bar {\Phi }}^n(f_1...f_k)(x^{(n)}) =
[\sum_{\alpha =0}^{k-1} {\hat {\pi }}^{\otimes \alpha }\otimes {\bar
{\Phi }} \otimes {\hat P}^{\otimes (k-\alpha -1)}]^n.
(f_1\otimes ...\otimes f_k)(x^{(n)})$ \\
and ${\bar {\Phi }}^n(f_1...f_k)\in C^0(U^{(n)},Y)$, where ${\hat
{\pi }}^{\otimes \alpha }\otimes {\bar {\Phi }}\otimes {\hat
P}^{\otimes (k-\alpha -1)}.(f_1\otimes ... \otimes f_k) := ({\hat
{\pi }}(f_1...f_{\alpha }))({\bar {\Phi }}f_{\alpha +1})({\hat
P}(f_{\alpha +2}...f_k))$ (see Lemma 7).}
\par {\bf Proof.} Applying operator ${\hat {\psi }}_n$ from Note 5
to both sides of Equation $7(1)$ we get Formula $(1)$ of this
Corollary.
\par {\bf 9. Lemma.} {\it Let $u\in C^{[n]}({\bf K}^s,{\bf K}^m)$,
$u({\bf K}^s)\subset U$ and $f\in C^{[n]}(U,Y)$, where $U$ is an
open subset in ${\bf K}^m$, $s, m\in \bf N$, $Y$ is a $\bf K$-linear
space, then
\par $(1)$ $(f\circ
u)^{[n]}(x^{[n]})=
[\sum_{j_1=1}^m...\sum_{j_n=1}^{m(n)}(A_{j_n,v^{[n-1]},t_n}...
A_{j_1,v^{[0]},t_1} f\circ u) (\Upsilon ^1\circ
p_{j_n}{\hat S}_{j_{n-1}+1,v^{[n-2]}t_{n-1}}$ \\
$...{\hat S}_{j_1+1,v^{[0]}t_1}u^{n-1}) (P_n\Upsilon ^1\circ
p_{j_{n-1}}{\hat S}_{j_{n-2}+1,v^{[n-3]}t_{n-2}}... {\hat
S}_{j_1+1,v^{[0]}t_1} u^{n-2})...(P_n...P_2\Upsilon ^1\circ
p_{j_1}u)+  \sum_{j_1=1}^m...\sum_{j_{n-1}=1}^{m(n-1)} ({\hat {\pi
}}^1(A_{j_{n-1},v^{[n-2]},t_{n-1}}...A_{j_1,v^{[0]},t_1}f\circ u)
[\sum_{\alpha =0}^{n-2} {\hat {\pi }}^{\otimes \alpha }\otimes
{\Upsilon } \otimes {\hat P}^{\otimes (n-\alpha -2)}] ((\Upsilon
^1\circ p_{j_{n-1}}{\hat S}_{j_{n-2}+1,v^{[n-3]}t_{n-2}}...{\hat
S}_{j_1+1,v^{[0]}t_1}u^{n-2})
\otimes ... \otimes (P_{n-1}...P_2\Upsilon ^1\circ p_{j_1}u))$ \\
$+ [\sum_{\alpha =0}^{n-2} {\hat {\pi }}^{\otimes \alpha }\otimes
{\Upsilon } \otimes {\hat P}^{\otimes (n-\alpha -2)}]
(\sum_{j_1=1}^m...\sum_{j_{n-2}=1}^{m(n-2)}({\hat {\pi
}}^1(A_{j_{n-2}, v^{[n-3]},t_{n-2}}...A_{j_1,v^{[0]},t_1}f\circ
u))\otimes [\sum_{\alpha =0}^{n-3} {\hat {\pi }}^{\otimes \alpha
}\otimes {\Upsilon } \otimes {\hat P}^{\otimes (n-\alpha -3)}]
((\Upsilon ^1\circ p_{j_{n-2}}{\hat
S}_{j_{n-3}+1,v^{[n-4]}t_{n-3}}... {\hat
S}_{j_1+1,v^{[0]}t_1}u^{n-3})\otimes ...\otimes
(P_{n-2}...P_2\Upsilon ^1\circ p_{j_1}u))+...$
\\ $+ [\sum_{\alpha =0}^2 {\hat {\pi }}^{\otimes \alpha }\otimes
{\Upsilon } \otimes {\hat P}^{\otimes (2-\alpha )}]^{n-3} \{
\sum_{j_1=1}^m\sum_{j_2=1}^{m(2)} ({\hat {\pi }}^1
A_{j_2,v^{[1]},t_2} A_{j_1,v^{[0]},t_1}f\circ u) ({\Upsilon
^1}\otimes {\hat P}^1 + {\hat {\pi }}^1\otimes {\Upsilon ^1})
((\Upsilon ^1\circ p_{j_2}{\hat S}_{j_1+1,v^{[0]}t_1}u)\otimes
(P_2\Upsilon ^1\circ
p_{j_1}u)) \} $ \\
$+ ({\Upsilon }\otimes {\hat P} + {\hat {\pi }}\otimes {\Upsilon
})^{n-2} \{ \sum_{j_1=1}^m({\hat {\pi }}^1A_{j_1,v^{[0]},t_1} f\circ
u)\otimes (\Upsilon ^2\circ p_{j_1}u) \} ](x^{[n]})$ \\
and $f\circ u\in C^0(({\bf K}^s)^{[n]},Y)$, where $S_{j,\tau
}u(y):=(u_1(y),...,u_{j-1}(y),u_j(y+\tau _{(s)}),u_{j+1}(y+ \tau
_{(s)}),...,u_m(y+\tau _{(s)}))$, $u=(u_1,...,u_m)$, $u_j\in \bf K$
for each $j=1,...,m$, $y\in {\bf K}^s$, $\tau =(\tau _1,...,\tau
_k)\in {\bf K}^k$, $k\ge s$, $\tau _{(s)}:=(\tau _1,...,\tau _s)$,
$p_j(x):=x_j$, $x=(x_1,...,x_m)$, $x_j\in \bf K$ for each
$j=1,...,m$, ${\hat S}_{j+1,\tau }g(u(y),\beta ):=g(S_{j+1,\tau
}u(y),\beta )$, $y\in {\bf K}^s$, $\beta $ is some parameter,
$A_{j,v,t}:=({\hat S}_{j+1,vt} \otimes t\Upsilon ^1\circ
p_j)^*\Upsilon ^1_j$, where $\Upsilon ^1$ is taken for variables
$(x,v,t)$ or corresponding to them after actions of preceding
operations as $\Upsilon ^k$, $\Upsilon ^1_jf(x,v_j,t) :=
[f(x+e_jv_jt)-f(x)]/t$, $(B\otimes A)^*\Upsilon ^1f_i\circ u^i
(x,v,t) := \Upsilon ^1_jf_i(Bu^i,v,Au^i)$, $B: {\bf K}^{m(i)}\to
{\bf K}^{m(i)}$, $A: {\bf K}^{m(i)}\to \bf K$, $e_j=
(0,...,0,1,0,...,0)\in {\bf K}^{m(i)}$ with $1$ on $j$-th place;
$m(i)=m+i-1$, $j_i=1,...,m(i)$; $u^1:=u$, $u^2:=(u^1,t_1\Upsilon
^1\circ p_{j_1}u^1)$,...,$u^n=(u^{n-1},t_{n-1}\Upsilon ^1\circ
p_{j_{n-1}}u^{n-1})$, $A_{j_1,v^{[0]},t_1}f\circ u =: f_1\circ u^1$,
$A_{j_n,v^{[n-1]},t_n}f_{n-1}\circ u^{n-1} =: f_n\circ u^n$, ${\hat
S}_*\Upsilon ^1f(z) := \Upsilon ^1f({\hat S}_*z)$.}
\par {\bf Proof.} At first consider $n=1$, then $(f\circ
u)^{[1]}(t_0,v,t) = [f(u(t_0+vt))-f(u(t_0))]/t$, where $t_0\in {\bf
K}^s$, $t\in \bf K$, $v\in {\bf K}^s$. Though we consider here the
general case mention, that in the particular case $s=1$ one has
$t_0\in \bf K$, $v\in \bf K$. Then \par $(f\circ u)^{[1]}(t_0,v,t)=
[f(u(t_0+vt)) - f(u_1(t_0), u_2(t_0+vt),...,u_m(t_0+vt))]/t +
[f(u_1(t_0),u_2(t_0+vt),u_3(t_0+vt),...,u_m(t_0+vt))-
f(u_1(t_0),u_2(t_0),u_3(t_0+vt),...,u_m(t_0+vt))]/t+...+
[f(u_1(t_0),...,u_{m-1}(t_0),u_m(t_0+vt)) - f(u(t_0))]/t$,\\ where
$u=(u_1,...,u_m)$, $u_j\in \bf K$ for each $j=1,...,m$. Since
$u_j(t_0+vt)-u_j(t_0) = tu_j^{[1]}(t_0,v,t)$, hence \par $(f\circ
u)^{[1]}(t_0,v,t) = \Upsilon
^1f((u_1(t_0),u_2(t_0+vt),...,u_m(t_0+vt)),e_1,t\Upsilon
^1u_1(t_0,v,t))$ \\  $\Upsilon ^1u_1(t_0,v,t) + \Upsilon
^1f((u_1(t_0),u_2(t_0),u_3(t_0+vt),...,u_m(t_0+vt)),e_2,t\Upsilon
^1u_2(t_0,v,t))$\\  $\Upsilon ^1u_2(t_0,v,t) +...+ \Upsilon
^1f(u(t_0),e_m,t\Upsilon ^1u_m(t_0,v,t))\Upsilon ^1u_m(t_0,v,t)$, \\
since $u_j\in \bf K$ for each $j=1,...,m$ and $\bf K$ is the field,
where $e_j=(0,...,0,1,0,...,0)\in {\bf K}^m$ with $1$ on $j$-th
place for each $j=1,...,m$. With the help of shift operators it is
possible to write the latter formula shorter:
\par $(2)$  $\Upsilon ^1(f\circ u)(y,v,t) = \sum_{j=1}^m
{\hat S}_{j+1,vt} \Upsilon ^1f(u(y),e_j,t\Upsilon ^1\circ
p_ju(y,v,t)) (\Upsilon ^1\circ p_ju(y,v,t))$, \\
where $p_j(x):=x_j$, $x=(x_1,...,x_m)$, $x_j\in \bf K$ for each
$j=1,...,m$, ${\hat S}_{j+1,\tau }g(u(y),\beta ):=g(S_{j+1,\tau
}u(y),\beta )$, $y\in {\bf K}^s$, $\tau \in {\bf K}^k$, $k\ge s$,
$\beta $ is some parameter. Introduce operators $A_{j,v,t}:=({\hat
S}_{j+1,vt} \otimes t\Upsilon ^1\circ p_j)^*\Upsilon ^1_j$, where
$\Upsilon ^1$ is taken for variables $(y,v,t)$ or corresponding to
them after actions of preceding operators as $\Upsilon ^k$
remembering that $y^{[k]}, v^{[k]}\in ({\bf K}^s)^{[k]}$, $t\in \bf
K$, $v^{[k]}=(v^{[k]}_1,v^{[k]}_2,v^{[k]}_3)$ with $v^{[k]}_1,
v^{[k]}_2\in ({\bf K}^s)^{[k-1]}$, $v^{[k]}_3\in {\bf K}^k$ for each
$k\ge 1$, in particular, $v^{[0]}=v^{[0]}_1$ for $k=0$, $\Upsilon
^1_jf(x,v,t) := [f(x+e_j v_jt)-f(x)]/t$, $(B\otimes A)^*\Upsilon
^1f_i\circ u^i (y,v,t) := \Upsilon ^1_jf_i(Bu^i,v,Au^i)$, $B: {\bf
K}^{m(i)}\to {\bf K}^{m(i)}$, $A: {\bf K}^{m(i)}\to \bf K$. For
example, in the particular case of $s=1$ we have $v^{[k]}\in ({\bf
K})^{[k]}$. Therefore, in the general case Formula $(2)$ takes the
form:
\par $(3)$ $\Upsilon ^1f\circ u(y,v,t)=\sum_{j=1}^m(A_{j,v,t}f\circ
u)(\Upsilon ^1\circ p_ju)(y,v,t)$.
\par Take now $n=2$, then
\par $\Upsilon ^2f\circ u(y^{[2]})= \Upsilon ^1 \sum_{j=1}^m
[(A_{j,v,t}f\circ u)(\Upsilon ^1\circ p_ju)(y,v,t)](y^{[2]})$. \\
In the square brackets there is the product, hence from Formula
$4(1)$ and Lemma 3 we get:
\par $(4)$ $\Upsilon ^2f\circ u(y^{[2]})=
\sum_{j=1}^m[(\Upsilon ^1A_{j,v^{[0]},t}f\circ u)(P_2\Upsilon
^1\circ p_ju) + ({\hat {\pi }}^1A_{j,v^{[0]},t}f\circ u)(\Upsilon
^2\circ p_ju)](y^{[2]})$. \\
Then from Formula $(3)$ applied to terms $A_{j,v,t}f\circ u$ it
follows, that $\Upsilon ^1A_{j_1,v^{[0]},t_1}f\circ u (y^{[2]})=
\sum_{j_2=1}^{m(2)}(A_{j_2,v^{[1]},t_2}A_{j_1,v^{[0]},t_1}f\circ u)
(\Upsilon ^1\circ p_{j_2} S_{j_1+1,v^{[0]}t_1}u)(y^{[2]})$, where
$v^{[0]}=v$, $t_1=t$ (see also Lemma 4). Therefore,
\par $(5)$ $\Upsilon ^2f\circ u(y^{[2]})=
[\sum_{j_1=1}^m \sum_{j_2=1}^{m(2)} (A_{j_2,v^{[1]},t_2}
A_{j_1,v^{[0]},t_1} f\circ u) (\Upsilon ^1 \circ p_{j_2}{\hat
S}_{j_1+1,v^{[0]}t_1}u)(P_2\Upsilon ^1\circ p_{j_1}u)
+\sum_{j_1=1}^m({\hat {\pi }}^1A_{j_1,v^{[0]},t_1}f\circ u)
(\Upsilon ^2 \circ p_{j_1}u)](y^{[2]})$. \\
Then for $n=3$ applying Formulas $(3)$ and $7(1)$ to $(5)$ we get:
\par $(6)$ $\Upsilon ^3f\circ u(y^{[3]})=
[\sum_{j_1=1}^m\sum_{j_2=1}^{m(2)}\sum_{j_3=1}^{m(3)}
(A_{j_3,v^{[2]},t_3} A_{j_2,v^{[1]},t_2} A_{j_1,v^{[0]},t_1}f\circ
u)$\\ $(\Upsilon ^1\circ p_{j_3}{\hat S}_{j_2+1,v^{[1]}t_2} {\hat
S}_{j_1+1,v^{[0]}t_1}u^2) (P_2\Upsilon ^1\circ p_{j_2}{\hat
S}_{j_1+1,v^{[0]}t_1}u) (P_3P_2\Upsilon ^1\circ p_{j_1}u) +$ \\
$\sum_{j_1=1}^m\sum_{j_2=1}^{m(2)} [({\hat {\pi
}}^1(A_{j_2,v^{[1]},t_2} A_{j_1,v^{[0]},t_1}f\circ u)) (\Upsilon
^2\circ p_{j_2}{\hat S}_{j_1+1,v^{[0]}t_1}u) (P_3P_2\Upsilon ^1\circ
p_{j_1}u) +$\\  $({\hat {\pi }}^1\{ (A_{j_2,v^{[1]},t_2}
A_{j_1,v^{[0]},t_1} f\circ u) (\Upsilon ^1\circ p_{j_2}{\hat
S}_{j_1+1,v^{[0]}t_1}u) \} (\Upsilon ^1P_2\Upsilon ^1\circ
p_{j_1}u)] +$\\  $\sum_{j_1=1}^m\sum_{j_3=1}^{m(3)}
(A_{j_3,v^{[2]},t_3}{\hat {\pi }}^1A_{j_1,v^{[0]},t_1}f\circ u)
(\Upsilon ^1\circ p_{j_3}{\hat S}_{j_1+1,v^{[0]}t_1}u) (P_3\Upsilon
^2\circ p_{j_1}u)$ \\  $+ \sum_{j_1=1}^m ({\hat {\pi
}}^2A_{j_1,v^{[0]},t_1}f\circ u)
(\Upsilon ^3\circ p_{j_1}u)] (y^{[3]})$. \\
Thus Formula $(1)$ is proved for $n=1, 2, 3$. Suppose that it is
true for $k=1,...,n$ and prove it for $k=n+1$. Applying Formula
$7(1)$ to both sides of $(1)$ we get:
\par $(7)$ $\Upsilon ^{n+1}f\circ u (y^{[n+1]}) =
[\sum_{j_1=1}^m...\sum_{j_{n+1}=1}^{m(n+1)}
(A_{j_{n+1},v^{[n]},t_{n+1}}... A_{j_1,v^{[0]}, t_1} f\circ u)$ \\
$(\Upsilon ^1\circ p_{j_{n+1}}{\hat S}_{j_n+1,v^{[n-1]}t_n}...{\hat
S}_{j_1+1,v^{[0]}t_1}u^n) (P_{n+1}\Upsilon ^1\circ p_{j_n}{\hat
S}_{j_{n-1}+1,v^{[n-2]}t_{n-1}}... {\hat S}_{j_1+1,v^{[0]}t_1}
u^{n-1})...$\\  $(P_{n+1}...P_2\Upsilon ^1\circ p_{j_1}u) +
\sum_{j_1=1}^m...\sum_{j_n=1}^{m(n)} ({\hat {\pi
}}^1(A_{j_n,v^{[n-1]},t_n}...A_{j_1,v^{[0]},t_1}f\circ u)$\\
$\Upsilon ^1((\Upsilon ^1\circ p_{j_n}{\hat
S}_{j_{n-1}+1,v^{[n-2]}t_{n-1}} ...{\hat
S}_{j_1+1,v^{[0]}t_1}u^{n-1}) ... (P_n...P_2\Upsilon ^1\circ
p_{j_1}u))+$\\ $\Upsilon
^1(\sum_{j_1=1}^m...\sum_{j_{n-1}=1}^{m(n-1)}({\hat {\pi
}}^1(A_{j_{n-1}, v^{[n-2]},t_{n-1}}...A_{j_1,v^{[0]},t_1}f\circ
u))\Upsilon ^1 ((\Upsilon ^1\circ p_{j_{n-1}}{\hat
S}_{j_{n-2}+1,v^{[n-3]}t_{n-2}}$\\  $... {\hat
S}_{j_1+1,v^{[0]}t_1}u^{n-2})...(P_{n-1}...P_2\Upsilon ^1\circ
p_{j_1}u))+... + \Upsilon ^{n-2}\{ \sum_{j_1=1}^m
\sum_{j_2=1}^{m(2)}$\\  $({\hat {\pi }}^1 A_{j_2,v^{[1]},t_2}
A_{j_1,v^{[0]},t_1}f\circ u)\Upsilon ^1((\Upsilon ^1\circ
p_{j_2}{\hat S}_{j_1+1,v^{[0]}t_1}u)(P_2\Upsilon ^1\circ p_{j_1}u))
\} +$ \\  $\Upsilon ^{n-1} \{ \sum_{j_1=1}^m{\hat {\pi
}}^1A_{j_1,v^{[0]},t_1} f\circ u)(\Upsilon ^2\circ p_{j_1}u)
\} ](y^{[n+1]})= $ \\
$[\sum_{j_1=1}^m...\sum_{j_{n+1}=1}^{m(n+1)}
(A_{j_{n+1},v^{[n]},t_{n+1}}... A_{j_1,v^{[0]},t_1} f\circ u)
(\Upsilon ^1\circ
p_{j_{n+1}}{\hat S}_{j_n+1,v^{[n-1]}t_n}$ \\
$...{\hat S}_{j_1+1,v^{[0]}t_1}u^n) (P_{n+1}\Upsilon ^1\circ
p_{j_n}{\hat S}_{j_{n-1}+1,v^{[n-2]}t_{n-1}}... {\hat
S}_{j_1+1,v^{[0]}t_1} u^{n-1})...(P_{n+1}...P_2\Upsilon ^1\circ
p_{j_1}u)+ \sum_{j_1=1}^m...\sum_{j_n=1}^{m(n)} ({\hat {\pi
}}^1(A_{j_n,v^{[n-1]},t_n}...A_{j_1,v^{[0]},t_1}f\circ u)
[\sum_{\alpha =0}^{n-1} {\hat {\pi }}^{\otimes \alpha }\otimes
{\Upsilon } \otimes {\hat P}^{\otimes (n-\alpha -1)}] ((\Upsilon
^1\circ p_{j_n}{\hat S}_{j_{n-1}+1,v^{[n-2]}t_{n-1}}...{\hat
S}_{j_1+1,v^{[0]}t_1}u^{n-1})\otimes ... \otimes(P_n...P_2\Upsilon
^1\circ p_{j_1}u))$ \\  $+ [\sum_{\alpha =0}^{n-1} {\hat {\pi
}}^{\otimes \alpha }\otimes {\Upsilon } \otimes {\hat P}^{\otimes
(n-\alpha -1)}] (\sum_{j_1=1}^m...\sum_{j_{n-1}=1}^{m(n-1)} ({\hat
{\pi }}^1(A_{j_{n-1}, v^{[n-2]},t_{n-1}}...A_{j_1,v^{[0]},t_1}f\circ
u))\otimes [\sum_{\alpha =0}^{n-2} {\hat {\pi }}^{\otimes \alpha
}\otimes {\Upsilon } \otimes {\hat P}^{\otimes (n-\alpha -2)}]
((\Upsilon ^1\circ p_{j_{n-1}}{\hat
S}_{j_{n-2}+1,v^{[n-3]}t_{n-2}}... {\hat
S}_{j_1+1,v^{[0]}t_1}u^{n-2})\otimes ...\otimes
(P_{n-1}...P_2\Upsilon ^1\circ p_{j_1}u))$
\\ $+ [\sum_{\alpha =0}^2 {\hat {\pi }}^{\otimes \alpha }\otimes
{\Upsilon } \otimes {\hat P}^{\otimes (2-\alpha )}]^{n-2} \{
\sum_{j_1=1}^m\sum_{j_2=1}^{m(2)} ({\hat {\pi }}^1
A_{j_2,v^{[1]},t_2} A_{j_1,v^{[0]},t_1}f\circ u) ({\Upsilon
^1}\otimes {\hat P}^1 + {\hat {\pi }}^1\otimes {\Upsilon ^1})
((\Upsilon ^1\circ p_{j_2}{\hat S}_{j_1+1,v^{[0]}t_1}u)\otimes
(P_2\Upsilon ^1\circ p_{j_1}u)) \} $ \\
$+ ({\Upsilon }\otimes {\hat P} + {\hat {\pi }}\otimes {\Upsilon
})^{n-1} \{ \sum_{j_1=1}^m({\hat {\pi }}^1A_{j_1,v^{[0]},t_1} f\circ
u)\otimes (\Upsilon ^2\circ p_{j_1}u) \} ](y^{[n+1]}).$ \\
Mention that in general $(\Upsilon ^{n+1}f\circ u)(y^{[n+1]})$ may
depend nontrivially on all components of the vector $y^{[n+1]}$
through several terms in Formula $(7)$. Thus Formula $(1)$ of this
Lemma is proved by induction.
\par {\bf 10. Corollary.} {\it Let $u\in C^n({\bf K}^s,{\bf K}^m)$,
$u({\bf K}^s)\subset U$ and $f\in C^n(U,Y)$, where $U$ is an open
subset in ${\bf K}^m$, $s, m\in \bf N$, $Y$ is a $\bf K$-linear
space, then
\par $(1)$ ${\bar {\Phi }}^n(f\circ
u)(x^{(n)})= [\sum_{j_1=1}^m...\sum_{j_n=1}^{m(n)}
(B_{j_n,v^{(n-1)},t_n}... B_{j_1,v^{(0)},t_1} f\circ u)$\\ $({\bar
{\Phi }}^1\circ p_{j_n}{\hat S}_{j_{n-1}+1,v^{(n-2)}t_{n-1}}...{\hat
S}_{j_1+1,v^{(0)}t_1}u^{n-1}) (P_n{\bar {\Phi }}^1\circ
p_{j_{n-1}}{\hat S}_{j_{n-2}+1,v^{(n-3)}_0,t_{n-2}}... {\hat
S}_{j_1+1,v^{(0)}t_1}u^{n-2})$ \\  $...(P_n...P_2{\bar {\Phi
}}^1\circ p_{j_1}u) + \sum_{j_1=1}^m...\sum_{j_{n-1}=1}^{m(n-1)}
({\hat {\pi }}^1(B_{j_{n-1},v^{(n-2)},t_{n-1}}...B_{j_1,v^{(0)},t_1}
f\circ u) [\sum_{\alpha =0}^{n-2} {\hat {\pi }}^{\otimes \alpha
}\otimes {\bar {\Phi }} \otimes {\hat P}^{\otimes (n-\alpha -2)}]$
\\ $(({\bar {\Phi }} ^1\circ p_{j_{n-1}}{\hat
S}_{j_{n-2}+1,v^{(n-3)}t_{n-2}}...{\hat
S}_{j_1+1,v^{(0)}t_1}u^{n-2})
\otimes ... \otimes(P_{n-1}...P_2{\bar {\Phi }}^1\circ p_{j_1}u))$ \\
$+ [\sum_{\alpha =0}^{n-2} {\hat {\pi }}^{\otimes \alpha }\otimes
{\bar {\Phi }} \otimes {\hat P}^{\otimes (n-\alpha -2)}]
(\sum_{j_1=1}^m...\sum_{j_{n-2}=1}^{m(n-2)}({\hat {\pi
}}^1(B_{j_{n-2},v^{(n-3)}, t_{n-2}}...B_{j_1,v^{(0)},t_1}f\circ
u))\otimes [\sum_{\alpha =0}^{n-3} {\hat {\pi }}^{\otimes \alpha
}\otimes {\bar {\Phi }} \otimes {\hat P}^{\otimes (n-\alpha -3)}]
(({\bar {\Phi }} ^1\circ p_{j_{n-2}}{\hat
S}_{j_{n-3}+1,v^{(n-4)}t_{n-3}}... {\hat
S}_{j_1+1,v^{(0)}t_1}u^{n-3})\otimes ...\otimes (P_{n-2}...P_2{\bar
{\Phi }}^1\circ p_{j_1}u))+...$
\\ $+ [\sum_{\alpha =0}^2 {\hat {\pi }}^{\otimes \alpha }\otimes
{\bar {\Phi }} \otimes {\hat P}^{\otimes (2-\alpha )}]^{n-3} \{
\sum_{j_1=1}^m\sum_{j_2=1}^{m(2)} ({\hat {\pi }}^1
B_{j_2,v^{(1)},t_2} B_{j_1,v^{(0)},t_1}f\circ u) ({\bar {\Phi
}}^1\otimes {\hat P}^1 + {\hat {\pi }}^1\otimes {\bar {\Phi }}^1)
(({\bar {\Phi }}^1\circ p_{j_2}{\hat S}_{j_1+1,v^{(0)}t_1}u)\otimes
(P_2{\bar {\Phi }}^1\circ p_{j_1}u)) \} $ \\
$+ ({\bar {\Phi }}\otimes {\hat P} + {\hat {\pi }}\otimes {\bar
{\Phi }})^{n-2} \{ \sum_{j_1=1}^m({\hat {\pi }}^1B_{j_1,v^{(0)},t_1}
f\circ u)\otimes ({\bar {\Phi }}^2\circ p_{j_1}u) \} ](x^{(n)})$ \\
and $f\circ u\in C^0(({\bf K}^s)^{(n)},Y)$ (see notation of Lemma
9), where $B_{j,v,t}:=({\hat S}_{j+1,vt} \otimes t{\bar {\Phi
}}^1\circ p_j)^*{\bar {\Phi }}^1_j$, where ${\bar {\Phi }}^1$ is
taken for variables $(x,v,t)$ or corresponding to them after actions
of preceding operations as ${\bar {\Phi }}^k$, ${\bar {\Phi
}}^1_jf(x,v,t) := [f(x+e_jv_jt)-f(x)]/t$, $(B\otimes A)^*{\bar {\Phi
}}^1f_i\circ u^i (x,v,t) := {\bar {\Phi }}^1_jf_i(Bu^i,v,Au^i)$, $B:
{\bf K}^{m(i)}\to {\bf K}^{m(i)}$, $A: {\bf K}^{m(i)}\to \bf K$,
$m(i)=m+i-1$, $j_i=1,...,m(i)$, $u^1=u$, $u^2:=(u^1,t_1{\bar {\Phi
}}^1\circ p_{j_1}u^1)$, $u^n:=(u^{n-1},t_{n-1}{\bar {\Phi }}^1\circ
p_{j_{n-1}}u^{n-1})$, ${\hat S}_* {\bar {\Phi }}^1f(x):= {\bar {\Phi
}}^1f({\hat S}_*x)$.}
\par {\bf Proof.} The restriction of operators of Lemma 9
on $W^{(n)}$ from Note 5 gives Formula $(1)$ of this corollary,
where $v^{(k)}\in ({\bf K}^s)^k\times {\bf K}^k$.
\par {\bf 11. Lemma.} {\it If $a\ne 0$, $a\in \bf K$, $U$
is an open subset in $X$, where $X$ and $Y$ are topological vector
spaces over $\bf K$, $f\in C^1(U,Y)$, $T\in \bf K$, $T\ne 0$, then
\par $(1)$ $\Upsilon ^1f (x,av,t/a)=a\Upsilon ^1 f(x,v,t)$ and
\par $(2)$ $\Upsilon ^1f(x,v,at)=a^{-1}\Upsilon ^1f(x,av,t)$ and
\par $(3)$ $\Upsilon ^1f(x/T,v,t)= T^{-1}\Upsilon ^1f(x/T,v,t/T)$
for each $(x,v,t)\in U^{[1]}$ and $(x,v,at)\in U^{[1]}$ and
$(x/T,v,t)\in U^{[1]}$ respectively.}
\par {\bf Proof.} We have identities: $\Upsilon ^1f(x,av,t/a)=
[f(x+vta/a)-f(x)]/(t/a)=a[f(x+vt)-f(x)]/t=a\Upsilon ^1f(x,v,t)$,
$\Upsilon ^1f(x,v,at)=[f(x+vta)-f(x)]/(at)=a^{-1}\Upsilon
^1f(x,av,t)$, for $g(x):=f(x/T)$ there is the equality $\Upsilon ^1
g(x,v,t)=[g(x+vt)-g(x)]/t=[f((x+vt)/T)-f(x/T)]/t=
T^{-1}[f(x/T+vt/T)-f(x/T)]/(t/T)= T^{-1}\Upsilon ^1f(x/T,v,t/T)$.
\par {\bf 12. Lemma.} {\it Let $u: {\bf K}\to {\bf K}^b$ be a
polynomial function: \par $(1)$ $u=\sum_{n=0}^ma_nx^n$, \\
where $a_n\in {\bf K}^b$ are expansion coefficients, $x\in \bf K$,
$m\in \bf N$, then
\par $(2)$ $\Upsilon ^qu(x^{[q]})=
\sum_{n=1}^ma_n\sum_{k_1=1}^n {n\choose k_1} \{ [
\sum_{k_2=1}^{n-k_1} {{n-k_1}\choose
k_2}...\sum_{k_q=1}^{n-k_1-...-k_{q-1}} {{n-k_1-...-k_{q-1}}\choose
k_q} $ \\  $x^{n-k_1-...-k_q} (\mbox{
}_1v^{[q-1]}_1)^{k_q}t_q^{k_q-1} {\cal S}_{v^{[q-1]},t_q}(\mbox{
}_1v^{[q-2]}_1)^{k_{q-1}}(t_{q-1})^{k_{q-1}-1} ...{\cal
S}_{v^{[1]},t_2}(v^{[0]})^{k_1}(t_1)^{k_1-1}] $ \\  $+
[x^{n-k_1-...-k_{q-1}} \sum_{k_2=1}^{k_1} {{n-k_1}\choose k_2} ...
\sum_{k_{q-1}=1}^{n-k_1-...-k_{q-2}} {{n-k_1-...-k_{q-2}}\choose
k_{q-1}} \sum_{k_q=1}^{k_{q-1}}{k_{q-1}\choose k_q}$\\ $ (\mbox{
}_1v^{[q-2]}_1)^{k_{q-1}-k_q} (\mbox{
}_1v^{[q-1]}_2)^{k_q}t_q^{k_q-1} {\cal
S}_{v^{[q-1]},t_q}(t_{q-1})^{k_{q-1}-1}$ \\  ${\cal
S}_{v^{[q-2]},t_{q-1}}(\mbox{
}_1v^{[q-3]}_1)^{k_{q-2}}(t_{q-2})^{k_{q-2}-1}... {\cal
S}_{v^{[1]},t_2}(v^{[0]})^{k_1}(t_1)^{k_1-1}]+...$ \\ $+ [
x^{n-k_1}(v^{[0]})^{k_1}\sum_{k_2=1}^{k_1-1} {{k_1-1}\choose k_2}...
\sum_{k_q=1}^{k_{q-1}-1} {{k_{q-1}-1}\choose k_q} t_q^{k_q-1}$\\
$(v_3^{[q-1]})^{k_q}t_{q-1}^{k_{q-1}-k_q-1}...
(v_3^{[2]})^{k_3}t_2^{k_2-k_3-1}(v_3^{[1]})^{k_2}t_1^{k_1-k_2-1}] \} $, \\
where $ {\cal S}_{v^{[q-1]},t_q}\mbox{ }_jx^{[q-1]}:=\mbox{
}_jx^{[q-1]}+\mbox{ }_jv^{[q-1]}t_q$ for each $j$, where
$x^{[q]}=(\mbox{ }_1x^{[q]},\mbox{ }_2x^{[q]},...)$ and this shift
operator acts on all terms on the right of it in a product.}
\par {\bf Proof.} In view of Lemma 3
\par $(3)$ $\Upsilon ^1u(x^{[1]})=\sum_{n=0}^ma_n((x+v^{[0]}t_1)^n-x^n)/t_1=
\sum_{n=1}^ma_n \sum_{k_1=1}^n {n\choose k_1}
x^{n-k_1}(v^{[0]})^{k_1}t_1^{k_1-1}$, \\
where ${n\choose k}$ are binomial coefficients,
\par $\Upsilon ^2u(x^{[2]})=\sum_{n=1}^ma_n \sum_{k_1=1}^n {n\choose k_1}
((x+v^{[1]}_1t_2)^{n-k_1}(v^{[0]}+v^{[1]}_2t_2)^{k_1}
(t_1+v_3^{[1]}t_2)^{k_1-1}
-x^{n-k_1}(v^{[0]})^{k_1}t_1^{k_1-1})/t_2$ in accordance with the
notation of the proof of Lemma 4. Then
\par $(4)$ $\Upsilon ^2u(x^{[2]})= \sum_{n=1}^ma_n\sum_{k_1=1}^n {n\choose k_1}
\{ [ \sum_{k_2=1}^{n-k_1} {{n-k_1}\choose k_2}
x^{n-k_1-k_2}(v^{[1]}_1)^{k_2}t_2^{k_2-1}
(v^{[0]}+v^{[1]}_2t_2)^{k_1}(t_1+v_3^{[1]}t_2)^{k_1-1}] +
[x^{n-k_1}\sum_{k_2=1}^{k_1} {k_1\choose k_2}(v^{[0]})^{k_1-k_2}
(v^{[1]}_2)^{k_2}t_2^{k_2-1}(t_1+v_3^{[1]}t_2)^{k_1-1}] + [
x^{n-k_1}(v^{[0]})^{k_1}\sum_{k_2=1}^{k_1-1} {{k_1-1}\choose k_2}
t_1^{k_1-k_2-1} (v_3^{[1]})^{k_2}t_2^{k_2-1}] \} $. \\
Therefore, Formulas $(3,4)$ prove Formula $(2)$ for $n=1$ and $n=2$.
Let formula $(2)$ be true for $n=1,...,q$, prove it for $n=q+1$.
Applying to both sides of Equation $(2)$ operator $\Upsilon ^1$ with
the help of Formula $7(2)$ or $7(1)$ we get Formula $(2)$ for
$n=q+1$ also.
\par {\bf 13. Corollary.} {\it Let suppositions of Lemma 13 be
satisfied, then $|\Upsilon ^qu(x^{[q]})|\le \max_{n=0}^m|a_n|$ for
each $x^{[q]}\in {\bf K}^{[q]}$ with $|x^{[q]}|\le 1$.}
\par {\bf Proof.} The absolute value of each term on the right side
of Formula $12(2)$ in the curled brackets is not greater than one,
since binomial coefficients are integer numbers and their
non-archimedean absolute value is not greater than one and each
component of the vector $x^{[q]}_j\in \bf K$ has an absolute value
not greater than one. Applying the non-archimedean inequality
$|y+z|\le \max (|y|,|z|)$ for arbitrary $y, z\in \bf K$ we get the
statement of this corollary.
\par {\bf 14. Corollary.} {\it Let $u$ be a polynomial as in Lemma
13, then
\par $(1)$ ${\bar {\Phi }}^qu(x^{(q)}) = \sum_{n=q}^m
a_n\sum_{k_1=1}^n\sum_{k_2=1}^{n-k_1}...\sum_{k_q=1}^{n-k_1-...-k_{q-1}}
{n\choose {k_1}} {{n-k_1}\choose {k_2}}$ \\
$...{{n-k_1-...-k_{q-1}} \choose {k_q}}
v_1^{k_1}...v_q^{k_q}t_1^{k_1-1}...t_q^{k_q-1}x^{n-k_1-...-k_q}$.}
\par {\bf 15. Lemma.} {\it  Let $V_j\in \bf R$, $V_j>0$, for
each $j\in \bf N$ and $\lim_{j\to \infty }V_j=0$. Suppose also that
$g\in C^{\infty }({\bf K}^l,{\bf K})$, there exists $R>0$ such that
$g(x)=0$ for each $|x|>R$, moreover, \par $(1)$ $|\Upsilon
^jg(x^{[j]})|\le C^{j+1}V_j^{-j}$ \\ for each $j$ and $|x^{[j]}|\le
R$, where $C>0$ is a constant. Put \par $u(x)=(a+ \sum_{k_1,k_2=0}^m
\sum_{i_1, i_2=1}^l \mbox{ }_{i_1,i_2}b_{k_1,k_2}
x_{i_1}^{k_1}x_{i_2}^{k_2})g(x/T)$ \\
for each $x\in {\bf K}^l$, where $T\in \bf K$, $0<|T|\le 1$, $a,
\mbox{ }_{i_1,i_2}b_{k_1,k_2} \in \bf K$. Then there exists a
constant $C_1>0$ independent of $a,\mbox{ }_{i_1,i_2}b_{k_1,k_2}, j,
x$ and $T$ such that
\par $(2)$ $\| \Upsilon ^ju(x^{[j]})\|_{C^0(B({\bf K}^{[j]},0,R),
{\bf K})} $ \\ $\le (\max_{i_1,i_2,k_1,k_2} (|a|,|\mbox{
}_{i_1,i_2}b_{k_1,k_2}|))\max (1,R^m) |T|^{-j} C_1^{j+1}V_j^{-j}$.}
\par {\bf Proof.} Apply Lemmas 4 and 11. For this calculate
by induction $\Upsilon ^1(a+bx)(x,v,t)=bv$, $\Upsilon
^2(a+bx)(x^{[2]})=bv^{[1]}_2$,..., $\Upsilon ^j(a+bx)(x^{[j]})=b
v^{[j-1]}_2$ for each $j\ge 3$. Therefore, $\| \Upsilon ^j (a+bx)
\|_{C^0(B({\bf K}^{[j]},0,R),{\bf K})}\le \max (|a|,|b|)R$ for each
$j\ge 0$. In general apply Formula $12(2)$ and Corollary 13. Then by
induction from Formula $11(3)$ it follows, that $\| \Upsilon
^jg(x/T) \|_{C^0(B({\bf K}^{[j]},0,R),{\bf K})}= |T|^{-j} \|
\Upsilon ^jg(x)\|_{C^0(B({\bf K}^{[j]},0,R),{\bf K})}$ for each
$j\ge 0$, where $\Upsilon ^0g=g$. Therefore, from Formula $4(1)$ and
the ultrametric inequality we have
\par $\| \Upsilon ^ju(x^{[j]})\|_{C^0(B({\bf K}^{[j]},0,R),{\bf K})}
\le $ \\  $(\max_{i_1,i_2,k_1,k_2} (|a|,|\mbox{
}_{i_1,i_2}b_{k_1,k_2}|))\max (1,R^m) \max_{k=0}^j|T|^{-k} \|
\Upsilon ^kg(x)\|_{C^0(B({\bf K}^{[k]},0,R),{\bf K})}\le $ \\
$( \max_{i_1,i_2,k_1,k_2} (|a|,|\mbox{ }_{i_1,i_2}b_{k_1,k_2}|))\max
(1,R^m) \max_{k=0}^j|T|^{-k}C^{k+1}V_k^{-k}$, \\
since $g(x)=0$ for $|x|>R$ and choosing $C_1>0$ such that $\infty
>C_1\ge \sup_{j=0}^{\infty
}[\sup_{k=0}^j(C^{k+1}V_k^{-k}V_j^j|T|^{j-k})^{1/(j+1)}]$ we get the
statement of this Lemma.
\par {\bf 16. Lemma.} {\it If $U$ is an open subset in ${\bf K}^b$,
$f: U\to \bf K$ is a marked function, then a space $Y_n$ of
functions $ \{ \Upsilon ^nf(x^{[n]}): v^{[0]}, \mbox{
}_lv^{[k]}_j\in \{ 0, 1 \}; j=1, 2; l=1,2,...; k=1,...,n-1 \} $ is
finite dimensional over $\bf K$ whenever it exists such that
$dim_{\bf K} Y_n \le (2^{m(n-1)}-1) dim_{\bf K}Y_{n-1}$, $n\in \bf
N$, $m(n)=2m(n-1)+1$ for $n\in \bf N$, $m(0)=b$.}
\par {\bf Proof.} We have the recurrence relation for a number of
variables belonging to $\bf K$, $m(n) = 2m(n-1)+1$ for each $n\in
\bf N$ corresponds to $\Upsilon ^nf(x^{[n]})$, $m(0)=b$ corresponds
to $f(x)$. For $n=1$ we have \par $\Upsilon
^1f(x,v,t)=(f(x+vt)-f(x))/t= [f(x+vt)-f(x+(v-\mbox{ }_bve_b)t)]/t +
[f(x+(v-\mbox{ }_bve_b)t)- f(x+(v-\mbox{ }_be_b-\mbox{
}_{b-1}e_{b-1})t)]/t+...+ [f(x+\mbox{ }_1ve_1t)-f(x)]/t$, where
$v=(\mbox{ }_1v,...,\mbox{ }_bv)$, $\mbox{ }_lv\in \bf K$ for each
$l=1,...,b$. We have that $\mbox{ }_lv\in \{ 0, 1 \} $ may take only
two values and the amount of such nonzero vectors $v$ is equal to
$2^b-1$. Thus the family $\{ \Upsilon ^1f(x+(v-\mbox{ }_bve_b- ...-
\mbox{ }_kve_k),\mbox{ }_kve_k,t): \mbox{ }_lv\in \{ 0, 1 \},
l=1,...,b \} $ of functions by $(x,t)\in {\bf K}^{b+1}$ spans over
$\bf K$ the space $\{ \Upsilon ^1f(x,v,t): \mbox{ }_lv\in \{ 0, 1
\}, l=1,...,b \} $. Its dimension over $\bf K$ for a given $f$ is
not greater, than $2^b-1$.
\par  Let the statement of this lemma be true for $n-1\ge 1$.
Then apply the operator $\Upsilon ^1$ to $\Upsilon
^{n-1}f(x^{[n-1]})$. Replacing in the proof above $f$ on $f^{[n-1]}$
we get the statement of this lemma for $n$ also, since $\Upsilon
^nf(x^{[n]})=\Upsilon ^1(\Upsilon
^{n-1}f(x^{[n-1]}))((x^{[n-1]})^{[1]})$ and
$x^{[n]}=(x^{[n-1]})^{[1]}$ considering $f^{[n]}(x^{[n]})$ by free
variables $(x,t_1,...,t_n)$.
\par {\bf 17. Corollary.} {\it For each $n\in \bf N$ and each
$b\in \bf N$ and a marked function $f: U\to \bf K$, where $U$ is
open in ${\bf K}^b$ there exists a finite system $\Lambda _n$ of
vectors $0\ne (y,v)$, $y, v\in {\bf K}^{m(n-1)}$ such that
\par $(1)$ $\sum_{(y,v)\in
\Lambda _n} C_{(y,v)} \Upsilon ^nf(x^{[n-1]}+y,
v,t_n)=0$ \\
is identically equal to zero as the function of $(x,t_1,...,t_n)$,
where $(x^{[n-1]}+y,v,t_n)\in U^{[n]}$, $x^{[n-1]} \in U^{[n-1]}$,
$v^{[0]}$ and $\mbox{ }_lv^{[k]}_j\in \{ 0, 1 \} $ for each $j, k,
l$, $y$ may depend on the parameters $t_1,...,t_n$ polynomially,
$0\ne C_{(y,v)}\in \bf K$ are constants for each $(y,v)$.}
\par {\bf Proof.} Take $card (\Lambda _n)>dim_{\bf K}Y_n$ and
$\mbox{ }_lv^{[k]}_j\in \{ 0, 1 \} $ for each $l=1,...,m(k-1)$,
$j=1, 2$ and $k=0,..., n-1$ such that $(x,v,t_1)\in U^{[1]}$. Then
\par $\Upsilon ^1f(x,v,t_1)=\Upsilon ^1f(x+(v-\mbox{
}_bve_b)t_1,\mbox{ }_bve_b,t_1)+\Upsilon ^1f(x+ (v-\mbox{
}_bve_b-\mbox{ }_{b-1}ve_{b-1})t,\mbox{
}_{b-1}ve_{b-1},t_1)+...+\Upsilon ^1f(x,\mbox{ }_1ve_1,t_1)$,\\
hence $\Upsilon ^f$ on vectors $\{ (x,v,t_1); (x+(v-\mbox{
}_bve_b)t_1,\mbox{ }_bve_b,t_1); (x+(v-\mbox{ }_bve_b-\mbox{
}_{b-1}ve_{b-1})t_1,\mbox{ }_{b-1}ve_{b-1},t_1);...; (x,\mbox{
}_1ve_1,t_1) \} $ is $\bf K$-linearly dependent system of functions
by $(x,t_1)$, where $\mbox{ }_lv\in \{ 0, 1 \} $, $l=1,...,b$,
$C_{(y,v)}\ne 0$. Let the statement be proved for $n-1$, then prove
it for $n$. Apply to both sides of equation \par
$\sum_{(y^{n-1},v^{n-1})\in \Lambda _{n-1}} C_{(y^{n-1},v^{n-1})}
\Upsilon
^{n-1}f(x^{[n-2]}+y^{n-1},v^{n-1},t_{n-1})=0$ \\
operator $\Upsilon ^1$, which is $\bf K$-linear, consequently, \par
$\sum_{(y^1,v^1)\in \Lambda _1} C_{(y^1,v^1)} \Upsilon
^1(\sum_{(y^{n-1},v^{n-1})\in \Lambda _{n-1}} \Upsilon
^{n-1}f)((x^{[n-2]}+y^{n-1},v^{n-1},t_{n-1})
+y^1,v^1,t_n) =0$, \\
where $\Lambda _1$, $y^1$ and $v^1$ already correspond to $\Upsilon
^{n-1}f(x^{[n-1]})$ instead of $f(x)$, we get Formula $(1)$ with
$C_{(y^n,v^n)}= C_{(y^1,v^1)}C_{(y^{n-1},v^{n-1})}\ne 0$ and with
$\Upsilon ^nf(x^{[n-2]}+y^{n-1},v^{n-1},t_{n-1})
+y^1,v^1,t_n)=\Upsilon ^nf(x^{[n-1]}+y^n,v^n,t_n)$.
\par {\bf 18. Corollary.} {\it If $U$ is an open subset in ${\bf K}^b$,
$f: U\to \bf K$ is a marked function, then a space $X_n$ of
functions $ \{ {\bar {\Phi }}^nf(x^{(n)}): \mbox{ }_lv_j\in \{ 0, 1
\}, l=1,...,b; j=1,...,n \} $ is finite dimensional over $\bf K$
whenever it exists such that $dim_{\bf K} X_n \le (2^b-1)^n$, $n\in
\bf N$. Moreover, there exists a finite system $\Lambda _n$ of
vectors $0\ne (y,v)$, $y \in {\bf K}^b$, $v\in ({\bf K}^b)^n$ such
that
\par $(1)$ $\sum_{(y,v)\in
\Lambda _n} C_{(y,v)} {\bar {\Phi }}^nf(x+y,
v,t_1,...,t_n)=0$ \\
is identically equal to zero as the function of $(x,t_1,...,t_n)$,
where $(x+y,v,t_1,...,t_n)\in U^{(n)}$, $x^{(n-1)} \in U^{(n-1)}$,
$v^{(0)}$ and $\mbox{ }_lv_j\in \{ 0, 1 \} $ for each $l=1,...,b$,
$j=1,...,n$, $y$ may depend on the parameters $t_1,...,t_n$
linearly, $0\ne C_{(y,v)}\in \bf K$ are constants for each $(y,v)$.}
\par {\bf Proof.} Restrict in the preceding formulas $\Upsilon
^nf(x^{[n]})$ on $W^{(n)}$ and from Lemma 25 and Corollary 26 we get
the statement of this corollary.
\par {\bf 19. Lemma.} {\it Let $U$ be an open subset in ${\bf K}^m$,
$Y$ be a $\bf K$-linear space. If $char ({\bf K})=0$, then either
$f\in C^{[n]}(U,Y)\cap C^{n+1}(U,Y)$ or $f\in C^{[n]}_b(U,Y)\cap
C^{n+1}_b(U,Y)$ if and only if either $f\in C^{[n+1]}(U,Y)$ or $f\in
C^{[n+1]}_b(U,Y)$. If $char ({\bf K})>0$, then $C^{[n]}(U,Y)\subset
C^n(U,Y)$ and $C^{[n]}_b(U,Y)\subset C^n_b(U,Y)$.}
\par {\bf Proof.} If $f\in C^{[n+1]}(U,Y)$ or $f\in C^{[n+1]}_b(U,Y)$,
then the restriction $\Upsilon ^{n+1}f|_{W^{(n+1)}}={\bar {\Phi
}}^{n+1}f$ is continuous or uniformly continuous on $V^{(n+1)}$
correspondingly, consequently, $f\in C^{n+1}(U,Y)$ or $f\in
C^{n+1}_b(U,Y)$ respectively. Since $C^{[n]}(U,Y)\subset
C^{[n+1]}(U,Y)$ or $C^{[n]}_b(U,Y)\subset C^{[n+1]}_b(U,Y)$, then
$f\in C^{[n]}(U,Y)\cap C^n(U,Y)$ or $f\in C^{[n]}_b(U,Y)\cap
C^n_b(U,Y)$ correspondingly (see also Note 5).
\par Let now $f\in C^{[n]}(U,Y)\cap
C^{n+1}(U,Y)$ or $f\in C^{[n]}_b(U,Y)\cap C^{n+1}_b(U,Y)$. For $n=0$
the statement of this lemma follows from Lemma 2. Suppose that the
statement of this lemma is true for $k=1,...,n$, then prove it for
$k=n+1$. In view of Lemma 9 we have, that $f^{[n+1]}(x^{[n]})$ has
the expression through the finite sum of terms $({\hat {\Phi
}}^{n+1}f\circ u^{n+1})h_{\beta }$ up to minor terms $({\Upsilon
}^if)h_{\beta }$ with $i\le n$, where $u^{n+1}\in C^{\infty }_b$ and
$h_{\beta }\in C^{\infty }_b$ are functions associated with $\bf
K$-linear and polynomial shift operators and their compositions
independent of $f$. We can write this in more details by induction.
Now consider the composite function: \par $(1)$ $g(x^{[q]}):=({\bar
{\Phi }}^q f)(u(x;\alpha ); e_{j_1},...,e_{j_q};a_1{\bar {\Phi
}}^{n(1)} u_{k_1}(x;\mbox{ }_1e;b_1),...,$\\  $a_s{\bar {\Phi
}}^{n(s)}u_{k_s}(x; \mbox{ }_se;b_s)) {\bar {\Phi
}}^{m(1)}w_1(x;\mbox{ }_1\xi
;c_1)... {\bar {\Phi }}^{m(r)}w_r(x;\mbox{ }_r\xi ;c_r)$ \\
appearing from the decomposition of $f^{[q]}$, where $a_1,...,a_s$
are polynomials of $t_1,...,t_q$ and $\mbox{ }_lv^{[k]}_j$,
$k=0,...,q$, $l=1,2,...$, $j=1, 2, 3$ for $k>0$, $b_l\subset \{
t_1,...,t_q \} $, $\mbox{ }_r\xi =(\mbox{ }_re_{i_1},...,\mbox{
}_re_{i_{m(r)}})$, $\mbox{ }_se= (\mbox{ }_se_{j_1},...,\mbox{
}_se_{j_{n(s)}})$; $k_j, s, r, n(s), m(r)\in \bf N$; $\alpha $ is a
parameter, $w_1,...,w_r$ are polynomials of $x, t_1,...,t_q, \mbox{
}_lv^{[k]}_j$; \\ here $\alpha , u, a_1,...,a_s, w_1,...,w_r, e_i,
\mbox{ }_j\xi , b_j$ are independent of $f$. The set of variables
$(x;t_1,...,t_q;\mbox{ }_lv^{[k]}_j: l=1,2,...; k=0,...,q; j=1, 2, 3
\} $ is in the bijective correspondence with the vector $x^{[q]}$.
Then act on this function $g$ by the operator $\Upsilon ^1$ at the
vector $x^{[q+1]}= (x^{[q]},v^{[q]},t_{q+1})$ such that $\Upsilon
^1g(x^{[q+1]})= [g(x^{[q]}+v^{[q]}t_{q+1})-g(x^{[q]})]/t_{q+1}$. For
the calculation of $\Upsilon ^1g$ apply Formulas $7(2)$ and $9(2)$
to $g$ and $({\bar {\Phi }}^q f)$ by all variables of functions in
this composition and product. It is nonlinear by ${\bar {\Phi
}}^qf$. As the result $\Upsilon ^1g$ is the $\bf K$-linear
combination of functions of the same type $(1)$ with $q+1$ instead
of $q$ and in general new functions in the composition and product
after actions on them operators $\Upsilon ^1$, $P_k$, ${\hat {\pi
}}$ and $S$. Shift operators $\hat S$ over $\bf K$ are infinite
differentiable and invertible such that for $\bf K$ with $char ({\bf
K})=0$ we have $\sum_{i=1}^k{\hat S}=k{\hat S}\ne 0$ and
$\sum_{i=1}^kt{\bar {\Phi }}^1id (x;v;t)=kt{\bar {\Phi }}^1id
(x;v;t)\ne 0$ for all $k\in {\bf N}= \{ 1, 2, 3,... \} $ and each
$t\ne 0$ and $v\ne 0$.
\par {\bf 20. Corollary.} {\it If $char ({\bf K})=0$, then a
function $f$ belongs to $C^{[n+1]}(U,Y)$ or $C^{[n+1]}_b(U,Y)$ if
and only if $f$ belongs to $C^{n+1}(U,Y)$ or $C^{n+1}_b(U,Y)$
respectively, moreover, there exists a constant $0<C_1<\infty $
independent of $f$ such that $ \| f\|_n\le \| f\|_{[n]}\le C_1\| f
\| _n$ for each $f\in C^{n+1}_b(U,Y)$, where $V^{[k]} := \{
x^{[k]}\in U^{[k]}: |v^{[q]}_1|=1; |\mbox{ }_lv^{[q]}_2t_{q+1}|\le
1, |v^{[3]}_q|\le 1 \quad \forall l, q \} $ or $V^{(k)} := \{
x^{(k)}\in U^{(k)}: |v_j|=1 \quad \forall j \} $ with norms either
\par $\| f\|_{[n]}:=\sup_{k=0,...,n; x^{[k]}\in V^{[k]}}
|f^{[k]}(x^{[k]})|$ or
\par $\| f\|_n := \sup_{k=0,...,n; x^{(k)}\in V^{(k)}}
|{\bar {\Phi }}^kf(x^{(k)})|$.}
\par {\bf Proof.} Apply Lemma 19 by induction for $k=1,...,n$
and use Lemma 2. If $g$ is a bounded continuous function $g: {\bf
K}^m\to \bf K$, then $[g(x+vt)-g(x)]/t$ is a bounded continuous
function by $(x,v,t)\in {\bf K}^m\times {\bf K}^m\times ({\bf
K}\setminus B({\bf K},0,\delta ))$, where $\delta >0$ is a constant.
If $L(X,Y)$ is the space of all bounded $\bf K$ linear operators $T:
X\to Y$ from a normed space $X$ into a normed space $Y$ over $\bf
K$, then operator norms $\| T \| _1:=\sup_{0\ne x\in X} \| Tx
\|_Y/\| x\|_X$, $\| T\|_2:= \sup_{0<|x|\le 1, x\in X}\| Tx \|_Y/\|
x\|_X$ and $\| T\|_3:= \sup_{|x|=1, x\in X}\| Tx \|_Y/\| x\|_X$ are
equivalent \cite{roo}. In view of Lemma 2 each operator ${\bar {\Phi
}}^jf(x;v_1,...,v_j;0,...,0)$ is $j$ multi-linear over $\bf K$ by
vectors $v_1,...,v_j\in X$. Therefore, the definition of the $C^n$
norm given above is worthwhile.  If $x^{[k]}\in V^{[k]}$, then
$|v^{[q]}_1|=1$ and $|\mbox{ }_lv^{[q-1]}_1+\mbox{
}_lv^{[q]}_2t_{q+1}|\le 1$ for each $l, q$. The inequality $\|
f\|_n\le \| f\|_{[n]}$ follows from ${\bar {\Phi
}}^nf=f^{[n]}_{W^{(n)}}$. The second inequality $\| f\|_{[n]}\le
C_1\| f \| _n$ follows from the decomposition of $f^{[q]}$ as a
finite $\bf K$-linear combination of terms having the form $19(1)$
for each $q=1,...,n$ and since norms of all terms are bounded and
expansion coefficients are independent of $f$, where $f^{[0]}=f$.
\par {\bf 21. Lemma.} {\it Let $U$ be an open subset in ${\bf K}^b$,
$b\in \bf N$, let also $f: U\to Y$ be a function with values in a
topological vector space $Y$ over $\bf K$. Then $f\in C^{[n]}(U,Y)$
or $f\in C^{[n]}_b(U,Y)$ or $f\in C^n(U,Y)$ or $f\in C^n_b(U,Y)$ if
and only if $\Upsilon ^kf(x^{[k]})\in
C^0(U^{[k]}_{j(0),j(1),...,j(k)},Y)$ or $\Upsilon ^kf(x^{[k]})\in
C^0_b(V^{[k]}_{j(0),j(1),...,j(k)},Y)$ or ${\bar {\Phi
}}^kf(x;e_{j(1)},...,e_{j(k)}; t_1,...,t_k)\in
C^0(U^{(k)}_{j(1),...,j(k)},Y)$ or ${\bar {\Phi
}}^kf(x;e_{j(1)},...,e_{j(k)}; t_1,...,t_k)\in
C^0_b(V^{(k)}_{j(1),...,j(k)},Y)$ for each $k=0,1,...,n$ and for
each $j(i)\in \{ 1,...,m_i \} $, $v^{[i]}=e_{j(i)}\in ({\bf
K}^b)^{[i]}$, $m_i=dim_{\bf K}({\bf K}^b)^{[i]}$, $i=0,1,...,k$,
$x^{[i+1]}= (x^{[i]},v^{[i]},t_{i+1})$ or respectively each
$j(1),...,j(k)\in \{ 1, 2,..., b \} $ with
$e_j=(0,...,0,1,0,...,0)\in {\bf K}^b$ is the vector with $1$ on the
$j$-th place, where $U^{[k]}_{j(0),...,j(l)}:=\{ x^{[k]}\in U^{[k]}:
v^{[i]}=e_{j(i)}, i=0,...,l \} $,
$V^{[k]}_{j(0),...,j(l)}=V^{[k]}\cap U^{[k]}_{j(0),...,j(l)}$,
$U^{(k)}_{j(1),...,j(l)} := \{ x^{(k)}\in U^{(k)}:
v_1=e_{j(1)},...,v_l=e_{j(l)} \} $, $V^{(k)}_{j(1),...,j(l)}=
V^{(k)}\cap U^{(k)}_{j(1),...,j(l)}$. Moreover, if each ${\Upsilon
}^kf(z^{[k]})|_{V^{[k]}_{j(0),...,j(k)}}$ or ${\bar {\Phi
}}^kf(z;e_{j(1)},...,e_{j(n)};t_1,...,t_n)$ is locally bounded, then
${\Upsilon }^kf(z^{[k]})$ or  ${\bar {\Phi }}^kf(z^{(k)})$ is
locally bounded respectively.}
\par {\bf Proof.} If $n=1$, then \par $(1)$ ${\bar {\Phi }}^1f(x;v_1;t_1)=
{\bar {\Phi }}^1f(\mbox{ }_1x,\mbox{ }_2x+\mbox{ }_2v_1t_1,...,
\mbox{ }_bx+\mbox{ }_bv_1t_1; e_1; \mbox{ }_1v_1t_1) \mbox{ }_1v_1 +
{\bar {\Phi }}^1f(\mbox{ }_1x,\mbox{ }_2x,\mbox{ }_3x+\mbox{
}_3v_1t_1,..., \mbox{ }_bx+\mbox{ }_bv_1t_1; e_2; \mbox{ }_2v_1t_1)
\mbox{ }_2v_1 +...$ \\  $+ {\bar {\Phi }}^1f(\mbox{
}_1x,...,\mbox{ }_bx; e_b; \mbox{ }_bv_1t_1) \mbox{ }_bv_1$, \\
hence ${\bar {\Phi }}^1f(x;v_1;t_1)\in C^0(U^{(1)},Y)$ or ${\bar
{\Phi }}^1f(x;v_1;t_1)\in C^0_b(V^{(1)},Y)$ if and only if ${\bar
{\Phi }}^1f(x;e_{j(1)};t_1)\in C^0(U^{(1)}_{j(1)},Y)$ or ${\bar
{\Phi }}^1f(x;e_{j(1)};t_1)\in C^0_b(V^{(1)}_{j(1)},Y)$ for each
$j(1)\in \{ 1,...,b \} $, where $x=(\mbox{ }_1x,...,\mbox{ }_bx)$,
$\mbox{ }_jx\in \bf K$ for each $j$, $\Upsilon ^1f={\bar {\Phi
}}^1f$. In accordance with Formula $10(1)$ or $9(1)$ we have the
expression of ${\bar {\Phi }}^kf(x;v_1,...,v_k;t_1,...,t_k)$ or
$\Upsilon ^kf(x^{[k]})$ throughout the sum of terms containing
${\bar {\Phi }}^kf(x;e_{j(1)},...,e_{j(k)};t_1,...,t_k)$ or
${\Upsilon }^kf(z^{[k]})|_{V^{[k]}_{j(0),...,j(k)}}$ with
multipliers belonging to $C^{\infty }_b(U,Y)$ or $C^{[\infty
]}_b(U,Y)$ putting in Formula $10(1)$ or $9(1)$ $u=id: {\bf K}^b\to
{\bf K}^b$, $id(x)=x$ for each $x$, $s=m=b$. From this the second
assertion follows. \par Suppose that the first statement of this
lemma is proved for all $k=0,1,...,n-1$. Then apply the operator
${\bar {\Phi }}^1$ to each ${\bar {\Phi
}}^{n-1}f(x;e_{j(1)},...,e_{j(n-1)};t_1,...,t_{n-1})$ and in
accordance with Formula $(1)$ with ${\bar {\Phi }}^{n-1}f$ or
$\Upsilon ^1$ to each $\Upsilon ^{n-1}f(x^{[n-1]})$ with
$x^{[n-1]}\in U^{[n-1]}_{j(0),...,j(n-1)}$ instead of $f$ we get the
same conclusion. Thus ${\bar {\Phi
}}^nf(x;e_{j(1)},...,e_{j(n-1)},v_n;t_1,...,t_{n-1},t_n)$ belongs to
$C^0$ or $C^0_b$ by its variables belonging to
$U^{(n)}_{j(1),...,j(n-1)}$ or to $V^{(n)}_{j(1),...,j(n-1)}$ or
$\Upsilon ^nf(x^{[n]})$ belongs to $C^0$ or $C^0_b$ by $x^{[n]}\in
U^{[n]}_{j(0),...,j(n-1)}$ or $x^{[n]}\in V^{[n]}_{j(0),...,j(n-1)}$
respectively if and only if ${\bar {\Phi
}}^{n-1}f(x;e_{j(1)},...,e_{j(n-1)},e_{j(n)};t_1,...,t_{n-1},t_n)$
belongs to $C^0(U^{(n)}_{j(1),...,j(n)},Y)$ or
$C^0_b(V^{(n)}_{j(1),...,j(n)},Y)$ or $\Upsilon
^nf(x^{[n]})|_{U^{[n]}_{j(0),...,j(n)}}\in C^{[n]}(
U^{[n]}_{j(0),...,j(n)},Y)$ or $\Upsilon
^nf(x^{[n]})|_{V^{[n]}_{j(0),...,j(n)}}\in C^{[n]}_b(
V^{[n]}_{j(0),...,j(n)},Y)$  respectively for each $j(n)$, where
$j(0),...,j(n)$ are arbitrary. Together with the induction
hypothesis this finishes the proof of this lemma.
\par {\bf 22. Lemma.} {\it Suppose that $U^k$ is open in
${\bf K}^k$ for each $2\le k\le m$ with $domain (f_k)=U^k$ and from
$f_k\circ u\in C^{[n]}({\bf K}^{k-1},Y)$ or $C^{[n]}_b({\bf
K}^{k-1},Y)$ or $C^n({\bf K}^{k-1},Y)$ or $C^n_b({\bf K}^{k-1},Y)$
for each $u\in C^{[\infty ]}({\bf K}^{k-1},{\bf K}^k)$ or
$C^{[\infty ]}_b({\bf K}^{k-1},{\bf K}^k)$ or $C^{\infty }({\bf
K}^{k-1},{\bf K}^k)$ or $C^{\infty }_b({\bf K}^{k-1},{\bf K}^k)$
with $image (u)\subset U^k$ it follows that $f_k\in C^{[n]}(U^k,Y)$
or $C^{[n]}_b(U^k,Y)$ or $C^n(U^k,Y)$ or $C^n_b(U^k,Y)$
respectively. Then for $domain (f)=U$ open in ${\bf K}^m$ from
$f\circ u\in C^{[n]}({\bf K},Y)$ or $C^{[n]}_b({\bf K},Y)$ or
$C^n({\bf K},Y)$ or $C^n_b({\bf K},Y)$ for each $u\in C^{[\infty
]}({\bf K},{\bf K}^m)$ or $C^{[\infty ]}_b({\bf K},{\bf K}^m)$ or
$C^{\infty }({\bf K},{\bf K}^m)$ or $C^{\infty }_b({\bf K},{\bf
K}^m)$ with $image (u)\subset U$ it follows, that $f\in C^{[n]}({\bf
K},Y)$ or $C^{[n]}_b({\bf K},Y)$ or $C^n({\bf K},Y)$ or $C^n_b({\bf
K},Y)$ respectively.}
\par {\bf Proof.} Write $f\circ u$ in the form $f\circ u=f\circ
u_{m-1}\circ u_{m-2}\circ ...\circ u_1$, where $u_j: {\bf K}^j\to
{\bf K}^{j+1}$ for each $j$ and $u: {\bf K}\to {\bf K}^m$ of
corresponding classes of smoothness. Applying supposition of this
lemma for $k=m, m-1,...,2$ we get that $f\circ u_{m-1}\circ ...\circ
u_j\in C^{[n]}({\bf K}^j,Y)$ provides $f\in C^{[n]}(U,Y)$ or also
for others classes of smoothness correspondingly for each
$j=m-1,m-2,...,1$.
\par {\bf 23. Lemma.} {\it Let $f: U\to {\bf K}^l$, where $U$ is open
in ${\bf K}^m$. Then $f\in C^{[n]}(U,{\bf K}^l)$ if and only if each
$\Upsilon ^nf(x^{[n]})_{U^{[n]}_{j(0),...,j(n)}}$ is continuous for
$v^{[k]}_3=0$ for each $k=1,...,n-1$.}
\par {\bf Proof.} In view of Lemma 21 it remains to prove, that
continuity of each $\Upsilon ^nf(x^{[n]})_{U^{[n]}_{j(0),...,j(n)}}$
is equivalent to the continuity of this family under the condition
$v^{[k]}_3=0$ for each $k=1,...,n-1$. Prove this by induction. We
already have that $v^{[k]}_3\in \{ 0, 1 \} $ for each $k=0,...,n-1$.
Denote by ${\hat S}_{n,t_n}$ the shift operator ${\hat
S}_{n,t_n}g(t_{n-1},\beta ):=g(t_{n-1}+t_n)$, where $\beta $ denotes
the family of all other variables of a function $g$. Then $\Upsilon
^nf(x^{[n]})=[{\hat S}_{n,t_n} \Upsilon
^{n-1}f(x^{[n-1]}+w^{[n-1]}t_n) - \Upsilon ^{n-1}f(x^{[n-1]})]/t_n =
[({\hat S}_{n,t_n}-I)\Upsilon ^{n-1}f(x^{[n-1]}+w^{[n-1]}t_n)]/t_n
+\Upsilon ^nf(x^{[n]})|_{v^{[n-1]}_3=0}$, \\
where $w^{[n-1]}$ differs from $v^{[n-1]}$ by $v^{[n-1]}_3$ such
that in $w^{[n-1]}$ it is zero and in $v^{[n-1]}$ it is one while
all others their components coincide such that
$x^{[n]}=(x^{[n-1]},v^{[n-1]},t_n)$, $x^{[n]}|_{v^{[n-1]}_3=0}=
(x^{[n-1]},w^{[n-1]},t_n)$. Since $[({\hat S}_{n,t_n}-I)\Upsilon
^{n-1}f(x^{[n-1]}+w^{[n-1]}t_n)]/t_n=[\Upsilon
^{n-1}f((x^{[n-2]},v^{[n-2]},t_{n-1}+t_n) + w^{[n-1]}t_n)- \Upsilon
^{n-1}f(x^{n-1}+w^{[n-1]}t_n)]/t_n$, where $x^{[0]}=x$,
$x^{[n-1]}=(x^{[n-2]},v^{[n-2]},t_{n-1})$, $n\ge 2$ and $k\ge 1$,
then $[({\hat S}_{n,t_n}-I)\Upsilon
^{n-1}f(x^{[n-1]}+w^{[n-1]}t_n)]/t_n=\Upsilon ^1_{s(n-1)}\Upsilon
^{n-1}f(x^{[n-1]}+v^{[n-1]}t_n) (t_{n-1}+t_n-t_{n-1})/t_n$ is
continuous, where $s(n-1)$ corresponds to the partial difference
quotients by the variable $t_{n-1}$. Then by induction get that
$({\hat S}_{k,t_k}-I)/t_k=\Upsilon ^1_{s(k)}$ for each $k=n-1,...,1$
which leads to the assertion of this lemma.
\par {\bf 24. Lemma.} {\it Suppose that $f\in C^n(U,Y)$ or
$f\in C^{[n]}(U,Y)$, where $U$ is open in ${\bf K}^m$, then each
${\bar {\Phi }}^nf(x^{(n)})$ has the symmetry by transposition of
pairs $(v_j,t_j)$ characterized by the Young tableaux consisting of
one row of length $n$, each $\Upsilon ^nf(x^{[n]})|_{\{ U^{[n]}:
v^{[k]}_3=0, k=1,...,n \} }$ is characterized by the Young tableaux
consisting of $2^{n-1}$ rows, where the first row of length $n$
contains numbers $1,...,n$, the second row of length $n-1$ contains
numbers $2,...,n$, the third and the fourth rows have lengths $n-2$
and contain numbers $3,...,n$ and so on, where the number of rows of
equal lengths $n-k$ is $2^{k-1}$ for $1\le k<n-1$. Moreover, if
$t_{i_1}=0$,...,$t_{i_l}=0$ as arguments of $\Upsilon ^nf$, then its
symmetry becomes higher with the amount of rows $2^{n-l}$ instead of
$2^{n-1}$.}
\par {\bf Proof.} The function ${\bar {\Phi
}}^nf(x;v_1,...,v_n;t_1,...,t_n)$ is symmetric relative to
transpositions $(v_i,t_i)\mapsto (v_j,t_j)$, since
$[(f(x+v_it_i+v_jt_j) - f(x+v_jt_j))/t_i
-(f(x+v_it_i)-f(x))/t_i]/t_j= [(f(x+v_it_i+v_jt_j) -
f(x+v_it_i))/t_j -(f(x+v_jt_j)-f(x))/t_j]/t_i$ for each $i\ne j$ and
so on by induction.
\par When $v^{[k]}_3=0$ for $1\le k\le n-1$ and $n\ge 2$ we have
\par $(1)$ $\Upsilon ^{k+1}f(x^{[k+1]})=
\{ [\Upsilon ^{k-1}f(x^{[k-1]}+(v^{[k-1]} + v^{[k]}_2t_{k+1})t_k +
v^{[k]}_1t_{k+1}) - \Upsilon
^{k-1}f(x^{[k-1]}+v^{[k]}_1t_{k+1})]/t_k - [\Upsilon
^{k-1}f(x^{[k-1]}+v^{[k-1]}t_k)-\Upsilon ^{k-1}f(x^{[k-1]})]/t_k \}
/t_{k+1} = \{ [\Upsilon ^{k-1}f(x^{[k-1]}+(v^{[k-1]} +
v^{[k]}_2t_{k+1})t_k+ v^{[k]}_1t_{k+1}) - \Upsilon ^{k-1}f(x^{[k-1]}
+ v^{[k-1]}t_k+ v^{[k]}_1t_{k+1})]/t_k + [\Upsilon ^{k-1}f(x^{[k-1]}
+ v^{[k-1]}t_k+ v^{[k]}_1t_{k+1})- \Upsilon
^{k-1}f(x^{[k-1]}+v^{[k]}_1t_{k+1})]/t_k - [\Upsilon
^{k-1}f(x^{[k-1]}+v^{[k-1]}t_k)-\Upsilon
^{k-1}f(x^{[k-1]})]/t_k \} /t_{k+1}$ \\
and this expression is symmetric relative to transpositions
$(v^{[k-1]},t_k)\mapsto (v^{[k]}_1,t_{k+1})$. Therefore, exclude
$v^{[k]}_3=0$ from the consideration such that $v^{[0]}:=v^{[0],1}$,
$v^{[1]}=(v^{[1],1},v^{[1],2},0)$, where $v^{[0],1}, v^{[1],1},
v^{[1],2}\in {\bf K}^m$. Then by induction define vectors
$v^{[k],i}\in {\bf K}^m$ such that $x^{[k]}+v^{[k]}t_{k+1}=
(x^{[k-1]}+v^{[k]}_1t_{k+1},v^{[k-1]}+v^{[k]}_2t_{k+1},
t_k+v^{[k]}_3t_{k+1})$ with $v^{[k]}_3=0$ and to this corresponds
$v^{[k-1],i}+v^{[k],i+2^{k-1}}t_{k+1}$ such that $v^{[k]}$ is
completely characterized  by $(v^{[k],i}: i=1,...,2^k)$, where $k\ge
1$.  Therefore, by induction $\Upsilon ^nf$ is symmetric relative to
transpositions $(v^{[k-1],i},t_k)\mapsto (v^{[k],i},t_{k+1})$ for
each $1\le i\le 2^{k-1}$, $1\le k\le n-1$. To $v^{[k],1}$ pose the
first row of length $n$ with numbers $1,...,n$ in boxes from left to
right, $k=0,1,...,n-1$. To vectors $v^{[k],i}$ with
$i=2^{k-1}+1,...,2^k$ and $k\ge 1$ pose rows in the Young tableaux
with such numbers in squares from left to right beginning with $k+1$
and ending with $n$ in each such $i$-th row.
\par If $t_{i_1}=0$,...,$t_{i_l}=0$ as arguments of $\Upsilon ^nf$,
then the symmetry of $\Upsilon ^nf$ up to notation corresponds to
$v^{[i_s-2]}+v^{[i_s-1]}_2t_{i_s}=v^{[i_s-2]}$ and $\Upsilon ^nf$ is
characterised by less amount of vectors $v^{[k],i}$, since $\Upsilon
^l_{t_{i_1},...,t_{i_l}}f={\bar {\Phi }}^l_{t_{i_1},...,t_{i_l}}f$
such that instead of $(v^{[k-1],j}: j=1,...,2^{k-1})$ it is
sufficient to take $j=1,...,2^{k-2}$ for $k=i_2$ for $k\ge 2$ and so
on excluding excessive vectors by induction on $s=3,...,l$.
\par {\bf 25. Lemma.} {\it Suppose that $f\in C^{n-1}(U,Y)$ or
$f\in C^{[n-1]}(U,Y)$, where $U$ is open in ${\bf K}^m$. Then
\par $(1)$ $f\in C^n(U,Y)$ or $f\in C^{[n]}(U,Y)$ if and only if ${\bar
{\Phi }}^nf(x;w,...,w;t_1,...,t_n)$ or \\ $\Upsilon ^nf(x^{[n]})|_{
\{ U^{[n]}: v^{[k],i}=w_s \quad \forall 2^{s-1}< i\le 2^s, 0\le s\le
k< n \} } $ is continuous for each marked $w\in {\bf K}^m$ or
$w_0,...,w_{n-1}\in {\bf K}^m$ respectively;
\par $(2)$ ${\bar {\Phi }}^nf$ or $\Upsilon ^nf$ is not locally
bounded if and only if there exists marked $w\in {\bf K}^m$ or
$w_0,...,w_{n-1}\in {\bf K}^m$ such that ${\bar {\Phi
}}^nf(x;w,...,w;t_1,...,t_n)$ or \\ $\Upsilon ^nf(x^{[n]})|_{ \{
U^{[n]}: v^{[k],i}=w_s \quad \forall 2^{s-1}< i\le 2^s, 0\le s\le k<
n \} } $ is not locally bounded.}
\par {\bf Proof.} In view of Lemma 11  and Formula $9(2)$ applied
by induction we have
\par ${\bar {\Phi }}^nf(x;w,...,w;t_1,...,t_n)=\sum_{i_1,...,i_n=1}^m
a_{i_1}...a_{i_n} {\bar {\Phi
}}^nf(x+t_1\sum_{l_1=i_1+1}^ma_{l_1}e_{i_1} +...
+t_n\sum_{l_n=i_n+1}^ma_{l_n}e_{l_n};e_{i_1},...,e_{i_n};a_{i_1}t_1,...,
a_{i_n}t_n)$ \\
for each $w=\sum_{i=1}^m a_ie_i$ if at least one $t_i\ne 0$, where
$a_i\in \bf K$, for convenience of notation
$\sum_{i=m+1}^ma_ie_i=0$. Then consider all $t_1,...,t_n\in \bf K$
such that $0\ne t_i\to 0$. Due to Lemma 24 and since $a_i$ are
arbitrary and can be taken nonzero, then each ${\bar {\Phi
}}^nf(x;e_{i_1},...,e_{i_n};a_{i_1}t_1,...,a_{i_n}t_n)$ is
continuous or locally bounded if and only if ${\bar {\Phi
}}^nf(x;w,...,w;t_1,...,t_n)$ is continuous or locally bounded for
each marked $w\in {\bf K}^m$. In view of Lemma 21 this provides
assertions $(1,2)$ for ${\bar {\Phi }}^nf$. \par We have ${\hat
P}^n(x^{[n]})|_{(U^{[n]}: v^{[k],i}=w_s \forall 2^{s-1}<i\le 2^s,
0\le s\le k<n )} = x+\sum_{k=0}^{n-1} \phi _{k+1}(t)w_k$, where
$\phi _l(t)=\sum_{1\le i_1<...<i_l\le n} t_{i_1}...t_{i_l}$ are
linearly independent symmetric polynomials, $l=1,...,n$,
$t=(t_1,...,t_n)$, in particular, $\phi _1(t)=t_1+...+t_n$. Put
$\alpha _{j,l}:=a_{j,s}$ for each $j=1,...,m$, $2^{s-1}<l\le 2^s$,
$s=0,...,n-1$, where $w_s=\sum_{i=1}^ma_{i,s}e_i$ with $a_{i,s}\in
\bf K$ for each $s=0,...,n-1$.
\par Applying Formula
$11(2)$ by induction we get
\par $\Upsilon ^nf(x^{[n]})|_{ \{ U^{[n]}: v^{[k],i}=w_s \quad \forall
2^{s-1}< i\le  2^s, 0\le s\le k< n \} } =
\sum_{i_0,...,i_{n-1}=1}^m\sum_{1\le q_k\le 2^k, k=0,...,n-1} $ \\ $
(\prod_{k=0}^{n-1}\alpha _{i_k,q_k}) \Upsilon ^nf(x^{[n]}_J)|_{ \{
U^{[n]}: v^{[s],l}=\delta _{l,q_s}e_{i_{s+1}}, \tau _{s+1}=\alpha
_{i_s,q_s}t_{s+1} \quad \forall s=0,...,n-1, 1\le l\le 2^s \} } $
\\ for each marked $w_s$ if at least one $t_i\ne 0$,
where $\delta _{i,j}=1$ for $i=j$ and $\delta _{i,j}=0$ for each
$i\ne j$, $J= \{ (i_k,q_k): k=0,...,n-1 \} $; ${\hat {\pi
}}^n(x^{[n]}_J)={\hat P}^n(y)$, where $\tau _{k+1}$ corresponds to
$x^{[n]}_J$ instead of $t_{k+1}$ for $x^{[n]}$, $y\in ({\bf
K}^m)^{[n]}$ corresponds to the set $(x;v^{[k],l}=\sum_{j_k\ge
i_k+\delta _{l,q_k}} \alpha _{j_k,l}e_{j_k}, k=0,...,n-1, 0\le s\le
k<n, 2^{s-1}<l\le 2^s ; t_1,...,t_n)$ in the notation introduced
above. Then consider all $t_1,...,t_n\in \bf K$ such that $0\ne
t_i\to 0$. Since $a_{i,s}\in \bf K$ are arbitrary constants which
can be taken nonzero, then from Lemmas 21 and 24 the statement of
this lemma for $\Upsilon ^nf$ as well follows, since $\Upsilon
^1(f^{[n-1]}(x^{[n-1]}))(x^{[n]})=
[f^{[n-1]}(x^{[n-1]}+v^{[n-1]}t_n)-f^{[n-1]}(x^{[n-1]})]/t_n$ and
$x^{[n]}=(x^{[n-1]},v^{[n-1]},t_n)$ and due to repeated application
of Formula $24(1)$, and $g(h(z)e_i,y)\in C^0$ by $(z,y)\in U_1\times
U_2$ is equivalent to $g(ue_i,y)\in C^0$ by $(u,y)\in h(U_1)\times
U_2$ for continuous function $h(z)$ by $z=(z_1,...,z_a)\in U_1$,
where $U_1$ and $U_2$ are domains in ${\bf K}^a$ and ${\bf K}^c$,
$g(ue_i,y)\in Y$, $h(U_1)\subset \bf K$.
\par {\bf 26. Lemma.} {\it If $f\in C^{n-1}(U,Y)$ or
$f\in C^{[n-1]}(U,Y)$, where $U$ is open in ${\bf K}^m$. Then ${\bar
{\Phi }}^nf(x^{(n)})|_{ \{ U^{(n)}: \exists i\quad |t_i|\ge \delta
\} }$ or $\Upsilon ^nf(x^{[n]})|_{ \{ U^{[n])}: \exists i\quad
|t_i|\ge \delta , v^{[k]}_3=0 \forall k \} }$ is continuous
respectively, where $\delta >0$.}
\par {\bf Proof.} Since ${\bar {\Phi }}^nf(x^{(n)})={\bar {\Phi }}^1(
{\bar {\Phi }}^{n-1}f(x^{(n-1)}))(x^{(n)})$ and $\Upsilon
^nf(x^{[n]})=\Upsilon ^1(\Upsilon ^{n-1}f(x^{[n-1]}))(x^{[n]})$
whenever it exists and $\Upsilon ^1f(x^{[1]})={\bar {\Phi
}}^1f(x^{(1)})=[f(x+vt)-f(x)]/t$ then in view of Lemmas 21 and 24 we
get the statement of this lemma, since ${\bar {\Phi
}}^{n-1}f(x^{(n-1)})$ or $\Upsilon ^{n-1}f(x^{[n-1]})$ is continuous
respectively and there is considered a domain with $|t_i|\ge \delta
$ and $t_i+v^{[i]}_3t_{i+1}=t_i$, where $v^{[i]}_3=0$.
\par {\bf 27. Lemma.} {\it Let $f: {\bf K}^b\to \bf K$ be a function
such that $f\circ u\in C^{[n]}({\bf K},{\bf K})$ or $f\circ u\in
C^n({\bf K},{\bf K})$ for $n\ge 0$ and $f\in C^{[n-1]}({\bf
K}^b,{\bf K})$ or $f\in C^{n-1}({\bf K}^b,{\bf K})$ for $n\ge 1$ for
each $u\in C^{[\infty ]}({\bf K},{\bf K}^b)$ or $u\in C^{\infty
}({\bf K},{\bf K}^b)$, where $\bf K$ is a field with a
non-archimedean valuation and $2\le b\in \bf N$, then ${\Upsilon
}^nf(x^{[n]})$ or ${\bar {\Phi }}^nf(x^{(n)})$ respectively is a
locally bounded function on $({\bf K}^b)^{[n]}$ or $({\bf
K}^b)^{(n)}$ and $f$ is continuous.}
\par {\bf Proof.} At first prove, that $f$ is continuous, when $n=0$,
since for $n\ge 1$ we have $C^0\subset C^{n-1}$. Suppose the
contrary, that there exists a sequence $\mbox{ }_jz$ such that
$\lim_{j\to \infty }\mbox{ }_jz=z_0$ and a limit of the sequence $\{
f(\mbox{ }_jz): j \} $ either does not exist or is not equal to
$f(z_0)$. Take $c_j$ and $r_j$ and $u(x)$ as above, then $\lim_{j\to
\infty }(f\circ u)(\mbox{ }_jx)=\lim_{j\to \infty } f(\mbox{
}_jz)\ne f(z_0)= (f\circ u)(y_0)$, hence $f\circ u$ is not
continuous at $y_0$ contradicting the assumption of this lemma.
\par Now suppose the contrary, that there exists $z_0^{[n]} \in ({\bf
K}^b)^{[n]}$ or $z_0^{(n)}\in ({\bf K}^b)^{(n)}$ such that $\Upsilon
^nf$ or ${\bar {\Phi }}^nf$ is unbounded in a neighborhood of
$z_0^{[n]}$ or $z_0^{(n)}$ correspondingly. As a neighborhood take a
ball $B(({\bf K}^b)^{[n]},z_0^{[n]},\epsilon )$ in $({\bf
K}^b)^{[n]}$ containing $z_0^{[n]}$ and of radius $\epsilon >0$ or
$B(({\bf K}^b)^{(n)},z_0^{(n)},\epsilon )$. Without loss of
generality we may suppose, that $z_0:=z_0^{[0]}=0\in {\bf K}^b$
making the shift $\phi (x):=f(x-z_0)$ when $z_0\ne 0$, where $z_0$
denotes the projection of $z_0^{[n]}$ in ${\bf K}^b$. Then there
exists a sequence $\mbox{ }_kz^{[n]}$ or $\mbox{ }_kz^{(n)}$ tending
to $z_0^{[n]}$ or $z_0^{(n)}$ when $k$ tends to the infinity such
that $\lim_{k\to \infty }|\Upsilon ^nf(\mbox{ }_kz^{[n]})|=\infty $
or $\lim_{k\to \infty }|{\bar {\Phi }}^nf(\mbox{ }_kz^{(n)})|=\infty
$ respectively, where $|x|=|x|_{\bf K}$ is the valuation in $\bf K$.
So we choose the sequence $\{ \mbox{ }_kz_0^{[n]}: k=1,2,... \} $
such that $|\mbox{ }_kv^{[n-1]}|\le 1$ and $|\mbox{ }_kt_j|\le 1$
for each $k\in \bf N$ and $j=1,...,n$. In view of Lemma 25 without
loss of generality there exists a marked $w\in {\bf K}^m$ or
$w_0,...,w_{n-1}\in {\bf K}^m$ such that ${\bar {\Phi
}}^nf(x;w,...,w;t_1,...,t_n)$ or \\ $\Upsilon ^nf(x^{[n]})|_{ \{
U^{[n]}: v^{[k],i}=w_s \quad \forall 2^{s-1}< i\le 2^s, 0\le s\le k<
n \} } $ is not locally bounded in a neighborhood of either
$z_0^{(n)}$ or $z_0^{[n]}$ with the sequence $\{ \mbox{ }_jz^{(n)}:
j\in {\bf N} \} $ or $\{ \mbox{ }_jz^{[n]}: j\in {\bf N} \} $ such
that either $\{ \mbox{ }_jz^{[n]}: j\in {\bf N}; \mbox{ }_jv_i=w,
i=1,...,n \} $ or $\{ \mbox{ }_jz^{[n]}: j\in {\bf N}; \mbox{
}_jv^{[k],i}=w_s \quad \forall 2^{s-1}< i\le 2^s, 0\le s\le k< n \}
$ respectively. At the same time due to Lemma 26 we can consider,
that $\lim_{j\to \infty } \max_{i=1}^n |\mbox{ }_jt_i| =0$. From
Formula $9(1)$ or $10(1)$ applied to $u=id$ and the conditions of
this lemma it follows, that all terms with orders $k<n$ of $B_*^kf$
or $A_*^k$ are continuous, hence there exists an ordered set $\{
j_n,...,j_1 \} $ such that the sequence either
\par $(1)$ $\{ (B_{j_n,v^{(n-1)},t_n}...
B_{j_1,v^{(0)},t_1} f\circ u) ({\bar {\Phi
}}^1\circ p_{j_n}{\hat S}_{j_{n-1}+1,v^{(n-2)}t_{n-1}}$ \\
$...{\hat S}_{j_1+1,v^{(0)}t_1}u^{n-1}) (P_n{\bar {\Phi }}^1\circ
p_{j_{n-1}}{\hat S}_{j_{n-2}+1,v^{(n-3)},t_{n-2}}... {\hat
S}_{j_1+1,v^{(0)}t_1} u^{n-2})...(P_n...P_2{\bar {\Phi }}^1\circ
p_{j_1}u) (\mbox{ }_jz_0^{(n)}): j\in {\bf N} \} $ or
\par $(2)$ $\{ (A_{j_n,v^{[n-1]},t_n}... A_{j_1,v^{[0]},t_1} f\circ u)
(\Upsilon ^1\circ p_{j_n}{\hat S}_{j_{n-1}+1,v^{[n-2]}t_{n-1}}$ \\
$...{\hat S}_{j_1+1,v^{[0]}t_1}u^{n-1}) (P_n\Upsilon ^1\circ
p_{j_{n-1}}{\hat S}_{j_{n-2}+1,v^{[n-3]}t_{n-2}}... {\hat
S}_{j_1+1,v^{[0]}t_1} u^{n-2})...(P_n...P_2\Upsilon ^1\circ
p_{j_1}u) (\mbox{ }_jz_0^{[n]}): j\in {\bf N} \} $ is unbounded for
$f=id$.
\par Now consider the same Formulas $10(1)$ or $9(1)$ for arbitrary $u$
satisfying conditions of this lemma. Again all terms with orders
$k<n$ of $B_*^kf\circ u$ or $A_*^kf\circ u$ are continuous and hence
bounded in a neighborhood of $z_0^{(n)}$ or $z_0^{[n]}$
respectively. We construct a curve $u$ in several steps leading to
the contradiction with the supposition of this lemma. \par  Mention
that $\Upsilon ^1 id(y,v^{[0]},t_1)=v^{[0]}$, where $y, v^{[0]}\in
{\bf K}^b$ and $t_1\in \bf K$. Then $\Upsilon ^2
id(y^{[2]})=(v^{[0]}+v^{[1]}_2t_2-v^{[0]})/t_2=v^{[1]}_2$ and
$\Upsilon ^3id_j(y^{[3]})=(\mbox{ }_jv^{[1]}_2+\mbox{
}_{j+b}v^{[2]}_2t_3-v^{[1]}_2)/t_3=\mbox{ }_{j+b}v^{[2]}_2$, where
$j=1,...,b$, $v^{[k]}=(\mbox{ }_1v^{[k]}_1,...,\mbox{ }_cv^{[k]}_1,
\mbox{ }_1v^{[k]}_2,...,\mbox{ }_cv^{[k]}_2,v_3^{[k]})$,
$c=c(k)=2^{k-1}-k+b(2^k-1)$, $\mbox{ }_jv^{[k]}_l\in \bf K$ for each
$j, k, l$, $id(y)=(id_1(y),...,id_b(y))=(y_1,...,y_b)$. Therefore,
we get by induction \par $\Upsilon ^mid_j(y^{[m]})=
\mbox{ }_{j(m)}v^{[m-1]}_2$, \\
for each $m\ge 2$, where $j(1)=j$, $j(2)=j$, $j(3)=j+b$,
$j(m)=j+2^{m-2}-(m-1)+b(2^{m-1}-1)$ for each $m\ge 4$, since
$j(m)=j+b+(2b+1)+ (2(2b+1)+1) + (2(2(2b+1)+1)+1) +...+
(2(2(...(2b+1)+1)+1)$ with $2$ in power $m-3$ in the latter term.
\par At first consider equations
\par $(3)$ ${\bar {\Phi }}^ku(x^{(n)})=\alpha _{k}{\bar
{\Phi }}^kid(z^{(n)})$  or
\par $(4)$ $\Upsilon ^ku(x^{[n]})=\alpha
_{k}\Upsilon ^kid(z^{[n]})$ for $k=0,1,...,n$ in neighborhoods of
$x_0^{(n)}$ and $z_0^{(n)}$ or $x_0^{[n]}$ and $z_0^{[n]}$ with
prescribed marked vectors $\eta $ or $\eta _0,...,\eta _{n-1}$ and
$w$ or  $w_0,...,w_{n-1}$ respectively, where $\eta $ or $\eta
_0,...,\eta _{n-1}$ are determined from the equations, $0\ne \alpha
_{k} \in \bf K$ are constants specified below for a sequence such
that $\lim_{j\to \infty }g_j=0$, where $0<q_j:=\min_{k=1}^n|\alpha
_{j,k}|\le g_j:=\max_{k=1}^n|\alpha _{j,k}|<1$. If $t_s=0$, then
equations for $\Upsilon ^ku$ simplify due to term $D_{t_s}$ instead
of $\Upsilon ^1_{t_s}$ for which $w_s$ does not play a role and we
can consider $\tau _s=0$, where $\tau _s$ play the same role for
$x^{(n)}$ and $x^{[n]}$ as $t_s$ for $z^{(n)}$ and $z^{[n]}$,
$s=1,...,n$. If $t_s\ne 0$, then we can take $\tau _s\ne 0$. In view
of Lemma 22 we can consider the data $(b-1,b)$ instead of $(1,b)$.
Since $w$ or $w_0,...,w_{n-1}$ are fixed vectors independent of $j$,
then we can resolve these equations for marked nonzero vectors $\eta
\in {\bf K}^{b-1}$ or $\eta _0,...,\eta _{n-1}\in {\bf K}^{b-1}$
corresponding to $\mbox{ }_jx^{(n)}$ or $\mbox{ }_jx^{[n]}$ such
that variables will be $\mbox{ }_jx\in {\bf K}^{b-1}$ and $\tau
_1,...,\tau _n$ for $u$ instead of $\mbox{ }_jz\in {\bf K}^b$ and
$t_1,...,t_n$ for $f$, such that $\lim_{j\to \infty
}\max_{i=1}^n|\tau _i|=0$. In view of Formulas $12(2)$ and $14(1)$
it is sufficient to consider a quadratic function
\par $(5)$ $u(h)=z + c \sum_{k_1,k_2=0}^2\sum_{i_1,i_2=1}^{b-1}\mbox{
}_{i_1,i_2}a_{k_1,k_2}h_{i_1}^{k_1}h_{i_2}^{k_2}$,\\ where $\mbox{
}_{i_1,i_2}a_{k_1,k_2}\in {\bf K}^b$, $c\in \bf K$,  $|\mbox{
}_{i_1,i_2}a_{k_1,k_2}|\le 1$ for each $i_1, i_2, k_1, k_2$,
$h=(h_1,...,h_{b-1})\in {\bf K}^{b-1}$. Thus we get
\par $(6)$ $|(B_{j_n,\eta ^{\otimes n},\tau _n}...
B_{j_1,\eta ,\tau _1} f\circ u) ({\bar {\Phi
}}^1\circ p_{j_n}{\hat S}_{j_{n-1}+1,\eta ^{\otimes (n-1)}\tau _{n-1}}$ \\
$...{\hat S}_{j_1+1,\eta \tau _1}u^{n-1}) (P_n{\bar {\Phi }}^1\circ
p_{j_{n-1}}{\hat S}_{j_{n-2}+1,\eta ^{\otimes (n-2)},\tau _{n-2}}...
{\hat S}_{j_1+1,\eta \tau _1} u^{n-2})$ \\  $...(P_n...P_2{\bar
{\Phi }}^1\circ p_{j_1}u) (\mbox{ }_jx_0^{(n)})|\ge |q_j|^n|\pi
|^{l_0+s_0} |(B_{j_n,w^{\otimes n},t_n}... B_{j_1,w,t_1} f\circ id)
({\bar {\Phi
}}^1\circ p_{j_n}{\hat S}_{j_{n-1}+1,w^{\otimes (n-1)}t_{n-1}}$ \\
$...{\hat S}_{j_1+1,wt_1}id^{n-1}) (P_n{\bar {\Phi }}^1\circ
p_{j_{n-1}}{\hat S}_{j_{n-2}+1,w^{\otimes (n-2)},t_{n-2}}... {\hat
S}_{j_1+1,w t_1} id^{n-2})$ \\  $...(P_n...P_2{\bar {\Phi }}^1\circ
p_{j_1}id) (\mbox{ }_jz_0^{(n)})|$ or
\par $(7)$ $|(A_{j_n,\eta ^{[n-1]},\tau _n}... A_{j_1,\eta ^{[0]},\tau _1}
f\circ u) (\Upsilon ^1\circ p_{j_n}{\hat S}_{j_{n-1}+1,\eta ^{[n-2]}
\tau _{n-1}}$ \\
$...{\hat S}_{j_1+1,\eta ^{[0]}\tau _1}u^{n-1}) (P_n\Upsilon ^1\circ
p_{j_{n-1}}{\hat S}_{j_{n-2}+1,\eta ^{[n-3]}\tau _{n-2}}... {\hat
S}_{j_1+1,\eta ^{[0]}\tau _1} u^{n-2})$ \\  $...(P_n...P_2\Upsilon
^1\circ p_{j_1}u) (\mbox{ }_jx_0^{[n]})| \ge |q_j|^n |\pi
|^{l_0+s_0}|(A_{j_n,w^{[n-1]},t_n}... A_{j_1,w^{[0]},t_1} f\circ id)
(\Upsilon ^1\circ p_{j_n}{\hat S}_{j_{n-1}+1,w^{[n-2]}t_{n-1}}$ \\
$...{\hat S}_{j_1+1,w^{[0]}t_1}id^{n-1}) (P_n\Upsilon ^1\circ
p_{j_{n-1}}{\hat S}_{j_{n-2}+1,w^{[n-3]}t_{n-2}}... {\hat
S}_{j_1+1,w^{[0]}t_1}id^{n-2})$ \\  $...(P_n...P_2\Upsilon ^1\circ
p_{j_1}id) (\mbox{ }_jz_0^{[n]})|$ \\
for each $j\in \bf N$, where $l_0\in \bf N$ is a marked number,
$s_0=s_0(j)\in \bf N$, each $w^{[k]}$ corresponds to marked
$w_0,...,w_{n-1}$, while $\eta , \eta _0, ...,\eta _{n-1}\in {\bf
K}^{b-1}$ are marked vectors for $u$, where $w^{\otimes
k}:=(w,...,w)\in X^{\otimes k}$ for $w\in X$ and $k\in \bf N$.
\par Take a function $\psi \in C^{\infty }({\bf K},{\bf K})$ such that
$\psi (x)=1$ for $|x|\le |\pi |$ and $\psi (x)=0$, when $|x|>|\pi
|$, for example, locally constant function, where $\pi \in \bf K$,
$0<|\pi |<1$. In particular, the characteristic function of $B({\bf
K},0,|\pi |)$ is locally analytic, since $\bf K$ is totally
disconnected with the base of its topology consisting of clopen
(closed and open simultaneously) balls, where $B(X,x,R) := \{ y\in
X, \rho (x,y)\le R \} $ for a topological space $X$ metrizable by a
metric $\rho $. It is proved further that such $\psi $ after
definite scalings suits construction below. Define now the functions
\par $u_j(h) :=(\xi _j\psi )((h-\mbox{ }_jx)/T_j)$, where
\par $\xi _j(h) := [\mbox{ }_{r_j}z_0 + c_j
\sum_{k_1,k_2=0}^2\sum_{i_1,i_2=1}^{b-1}\mbox{
}_{i_1,i_2}a_{k_1,k_2}h_{i_1}^{k_1}h_{i_2}^{k_2}]$ such that $\xi
_j(0)=\mbox{ }_{r_j}z_0$ and put
\par $u(x) := \sum_{j=1}^{\infty }u_j(x)$, \\ where $x=
(x_1,...,x_{b-1})\in {\bf K}^{b-1}$, each $\mbox{
}_{i_1,i_2}a_{k_1,k_2}\in {\bf K}^b$ is marked, $c_j\in {\bf
K}\setminus \{ 0 \} $. Choose $r_j\in {\bf N}$, $\mbox{ }_jx_i,
T_j\in \bf K$ later on. All $u_j$ have disjoint supports, hence the
series is convergent, $u$ is of class $C^{[\infty ]} :=
\bigcap_{n=1}^{\infty } C^{[n]}$ in ${\bf K}\setminus \{ z_0 \} $.
\par Consider the sets $\lambda _i := \{ j\in {\bf N}:
\mbox{ }_jt_i=0 \} $, then either $card ({\bf N}\setminus \lambda
_i) = \aleph _0$ or $card (\lambda _i)=\aleph _0$ or both
cardinalities are $\aleph _0$. Consider intersections $A_1\cap
...\cap A_n$, where $A_i=\lambda _i$ or $A_i={\bf N}\setminus
\lambda _i$. The union of all such finite intersections is $\bf N$.
Therefore, one of these intersections is of the cardinality $\aleph
_0$. Thus, there exists a subsequence $\{ j(l): l\in {\bf N} \} $
such that $\mbox{ }_st_{j(l)}=0$ for each $l$ and every $s\in \{
i_1,...,i_r \} $ and $\mbox{ }_st_{j(l)}\ne 0$ for each $s\in \{
1,...,n \} \setminus \{ i_1,...,i_r \} $, where $0\le r\le n$. After
the enumeration we can consider a sequence with such property. For
such a sequence we can choose a subsequence which after enumeration
has the property:
\par $(8)$ $|\mbox{ }_{j+1}t_i|\le |\pi |^{s(j)}|\mbox{ }_jt_i|$ and
\par $(9)$ $|\pi |^{r(j)}b_{j+1}\ge b_j$ for each $j\in \bf N$
and $i\in \{ 1,...,n \} $, \\ where $s(j), r(j)\in \bf N$ are
sequences specified below;
\par $b_j:=|(B_{j_n,w^{\otimes n},t_n}... B_{j_1,w,t_1} f\circ id)
({\bar {\Phi }}^1\circ p_{j_n}{\hat S}_{j_{n-1}+1,
w^{\otimes (n-1)}t_{n-1}}$ \\
$...{\hat S}_{j_1+1,w\tau _1}id^{n-1}) (P_n{\bar {\Phi }}^1\circ
p_{j_{n-1}}{\hat S}_{j_{n-2}+1,w^{\otimes (n-2)},t_{n-2}}... {\hat
S}_{j_1+1,w t_1}id^{n-2})...(P_n...P_2{\bar {\Phi }}^1\circ
p_{j_1}id) (\mbox{ }_jz_0^{(n)})|$ \\
or with analogous Properties $(8,9)$ for $A_*^nf$ instead of
$B_*^nf$.
\par Now choose $r_j$ and $c_j$ such that
$\lim_{j\to \infty } c_jT_j^{-q}=0$ for each $q\in \bf N$, for
example, $c_j=T_j^j$, where $\lim_{j\to \infty } T_j = 0 $,
$|T_j|>|T_{j+1}|$ for each $j$, $T_j\ne 0$ for each $j$. Then choose
$r_j\in \bf N$ such that $\max_{l=1}^n(|\mbox{ }_{r_j}t_l|)\le
|c_j|$ and $|\mbox{ }_{r_j}z_0|)\le |c_j|$ for each $j$ and
$\lim_{j\to \infty }|c_j^n \Upsilon ^nf(\mbox{
}_{r_j}z_0^{[n]})|=\infty $. Take $\mbox{ }_jx\in {\bf K}^{b-1}$
such that $\mbox{ }_jx_i=(\pi ^{-1}\sum_{k=1}^{j-1}T_k)+T_j$, where
$\pi \in \bf K$, $0<|\pi |<1$, $i=1,...,b-1$. Since
$|T_j|>|T_{j+1}|>0$ for each $j\in \bf N$, then $|\mbox{ }_jx-\mbox{
}_{j+1}x|=|T_{j+1}+T_j(\pi ^{-1}-1)|=|\pi ^{-1}T_j|>|T_j|$ and
$|\mbox{ }_kx-\mbox{ }_{k+1}x|\ge \min (|\mbox{ }_kx-\mbox{
}_{k+1}x|,|\mbox{ }_{k+1}x-\mbox{ }_{k+2}x|,...,|\mbox{ }_jx-\mbox{
}_{j+1}x|)\ge \min (|T_k|,...,|T_j|)$ for each $k\le j$,
consequently, $B({\bf K}^{b-1},\mbox{ }_kx,|T_k|)\cap B({\bf
K}^{b-1},\mbox{ }_{j+1}x,|T_{j+1}|)=\emptyset $ for each $k\le j$,
hence $supp (u_j)\cap supp (u_k)=\emptyset $ for each $k<j$. Take
$s_0(j+1)\ge s_0(j)+j+1$ and $|\pi |^{s_0(j+1)}<|q_j|\le g_j\le |\pi
|^{s_0(j)}$ and $r(j)\ge s_0(j)2n$ and $s(j)\ge s_0(j)$ for each
$j$, where $l_0$ is such that $0<|\pi |^{l_0}<1/2$ (see also
$(6-9)$).
\par Denote $y_0:=\lim_{j\to \infty }\mbox{
}_jx$. Then $u$ is of class $C^{[\infty ]}_b$ or $C^{\infty }_b$ in
a neighborhood of $z_0$. To prove this we show, that $\Upsilon
^qu(z^{[q]})$ or ${\bar {\Phi }}^qu(z^{(q)})$ tends to zero as $z$
tends to $z_0=0$, where $|z^{[q]}|<\epsilon $ or $|z^{(q)}|<\epsilon
$, since then for $ |t_1|\ge \epsilon ,$ or $...,|t_n|\ge \epsilon
$, $|v^{[q-1]}|\le 1$ the continuity will be evident. For this we
use Lemma 15. Mention, that $\| \Upsilon ^q\psi (x) \|_{C^0(B({\bf
K}^{[q]},0,R),{\bf K})}<\infty $ for each $q$ and each $R>0$.
Indeed, $\max_{x\in \bf K} |\psi (x)|=1$ so that $\Upsilon ^0\psi $
is bounded for $q=0$. For $q=1$ we have $\Upsilon ^1\psi
(x,v,t_1)=0$ for $\max (|x|,|x+vt_1|)\le |\pi |$ or $\min
(|x|,|x+vt_1|)>|\pi |$, $\Upsilon ^1\psi (x,v,t_1)=1/t_1$ for either
$|x|\le |\pi |$ and $|x+vt_1|>|\pi |$ or $|x|>|\pi |$ and
$|x+vt_1|\le |\pi |$. Since we consider the domain $|x^{[1]}|\le R$,
then $|v|\le R$, consequently, $\| \Upsilon ^1\psi (x)
\|_{C^0(B({\bf K}^{[1]},0,R),{\bf K})}\le R |\pi |^{-1}$, since
$|t_1|^{-1}\le |\pi |^{-1}R$ in the considered domain, when
$\Upsilon ^1\psi (x,v,t_1)\ne 0$. The function $\Upsilon ^1\psi
(x,v,t_1)$ is the product of the locally constant function by
variables $(x,v)$ and the function $1/t_1$ with $|\pi /v|\le
|t_1|\le R$, when this function is nonzero and $v\ne 0$, hence $|\pi
|/R\le |v|\le R$, that is $|\pi |/R\le |t_1|\le R$, where $\Upsilon
^1\psi (x,0,t_1)=0$ for each $x$ and $t_1$. Evidently, by induction
that $\Upsilon ^q\psi (x^{[q]})$ is in $C^0(B({\bf
K}^{[q]},0,R),{\bf K})$  with the finite norm $\| \Upsilon ^q\psi
(x^{[q]})\|_{C^0(B({\bf K}^{[q]},0,R),{\bf K})}\le
C^{q+1}V_q^{-q}\le (q+1)(R/|\pi |)^q$ with $V_q=1$ and
$C:=\lim_{q\to \infty }[(q+1)(R/|\pi |)^q]^{1/(q+1)}$ for non-scaled
$\psi $ for each $q\in \bf N$ and each $R\ge 1$. In general for
scaled $\psi $ put $V_q:=\min_{j=1}^q|T_j|>0$. At the same time for
each $x$ with $|x-\mbox{ }_jx|\le |T_j|$ and $|v^{[k]}|\le R$ and
$|t_{k+1}|\le R$ for each $k=0,...,n-1$ in accordance with Lemmas 4,
12, 15 and Corollary 13
\par $(10)$ $| \Upsilon
^qu(x^{[q]})| \le (\max (1,R^2))|c_j||T_j|^{-q}C_1^{q+1} V_q^{-q}$ \\
which tends to zero as $j$ tends to the infinity, since $C_1^{q+1}
\le (q+1)(R/|\pi |)^q$, $0<|T_{j+1}|<|T_j|$ for each $j$ and
$\lim_{j\to \infty } c_jT_j^{-\beta }=0$ for every $\beta \in \bf
N$, where $R\ge 1$.
\par If each term in Formula $9(1)$ OR $10(1)$ would be locally
bounded, then ${\bar {\Phi }}^n(f\circ u)(\mbox{ }_{r_j}x^{(n)})$ or
$\Upsilon ^n(f\circ u)(\mbox{ }_{r_j}x^{[n]})$ would be locally
bounded. Since each ${\bar {\Phi }}^kf$ or $\Upsilon ^kf$ is locally
bounded for $k<n$ by our supposition above, then from Formula $9(1)$
or $10(1)$ and the condition $\lim_{j\to \infty }|c_j^n \Upsilon
^nf(\mbox{ }_{r_j}z_0^{[n]})|=\infty $ or $\lim_{j\to \infty }|c_j^n
{\bar {\Phi }}^nf(\mbox{ }_{r_j}z_0^{(n)})|=\infty $ it follows,
that there exists a term or a finite sum of terms of the type \par
$(A_{j_n,v^{[n-1]},t_n}... A_{j_1,v^{[0]},t_1}f\circ u) (\Upsilon
^1\circ p_{j_n}
S_{j_{n-1}+1,v^{[n-2]}t_{n-1}}$ \\
$...S_{j_1+1,v^{[0]}t_1}u^{n-1}) (P_n\Upsilon ^1\circ
p_{j_{n-1}}S_{j_{n-2}+1,v^{[n-3]}t_{n-2}}... S_{j_1+1,v^{[0]}t_1}
u^{n-2})...(P_n...P_2\Upsilon ^1\circ
p_{j_1}u)\mbox{ }_{r_l}x^{[n]})$ \\
which absolute value tends to the infinity  for a particular set
$\omega $ of indices $(j_1,...,j_n)$ and a subsequence $\{ \mbox{
}_{r_l}x^{[n]} : j\in {\bf N} \} $ or analogously for $B_*^nf\circ
u$ instead of $A_*^nf\circ u$. But this contradicts supposition of
this lemma in view of Lemmas 9, 21 and Corollary 10. Therefore,
$\Upsilon ^nf$ or ${\bar {\Phi }}^nf$ respectively is locally
bounded.
\par {\bf 28. Remark.} Though $u\in C^{\infty }({\bf K},{\bf K}^b)$,
but $u$ is not locally analytic in general, since the sequence $\{
x_j: j \} $ converges to $y_0\in \bf K$ and $u$ has not a series
expansion in a neighborhood of $y_0$ with positive radius of
convergence.
\par {\bf 29. Definitions.} Let $\phi : (0,\infty )\to (0,\infty )$
be a function such that $\lim_{q\to 0}\phi (q)=0$. By either ${\bf
K}(\phi )$ or ${\bf K}(u,\phi )$ we denote the $\bf K$-linear space
of all functions $f: {\bf K}^m \to \bf K$ such that for each bounded
subset $U$ in ${\bf K}^m$ there exists a constant $C>0$ such that
either
\par $(1)$ $|f(x+y)-f(x)|
\le C\phi (|y|)$, \\
when $x\in U$ and $x+y\in U$ or
\par $(2)$ $|f(x+ut)-f(x)|
\le C\phi (|t|)$,\\
when $x\in U$ and $x+ut\in U$ respectively, where $u\in {\bf K}^m$
is a nonzero vector. In the particular case of $\phi (q)= q^w$,
where $0<w\le 1$, we also denote ${\bf K}(q^w)=: Lip (w)$ and ${\bf
K}(u,q^w)=: Lip (u,w)$.
\par Then we denote by $C^{[n],w}({\bf K}^m,{\bf K})$ or
$C^{[n]}_{\phi }({\bf K}^m,{\bf K})$ or $C^{n,w}({\bf K}^m,{\bf K})$
or $C^n_{\phi }({\bf K}^m,{\bf K})$ the $\bf K$-linear space of all
functions $f\in C^{[n]}({\bf K}^m,{\bf K})$ in the first and the
second cases or in $C^n({\bf K}^m,{\bf K})$ in the third and the
fourth cases such that $f^{[n]}(x^{[n]})\in Lip (w)$ or
$f^{[n]}(x^{[n]})\in {\bf K}(\phi )$ or ${\bar {\Phi
}}^nf(x^{(n)})\in Lip (w)$ or ${\bar {\Phi }}^nf(x^{(n)})\in {\bf
K}(\phi )$ respectively.
\par {\bf 30. Lemma.} {\it Let suppositions of Lemma 27 be satisfied
and moreover $\Upsilon ^n(f\circ u)\in {\bf K}(\phi )$ or ${\bar
{\Phi }}^n(f\circ u)\in {\bf K}(\phi )$ for each $u\in C^{[\infty
]}({\bf K},{\bf K}^b)$ or $C^{[\infty ]}_b({\bf K},{\bf K}^b)$ or
$u\in C^{\infty }({\bf K},{\bf K}^b)$ or $C^{\infty }_b({\bf K},{\bf
K}^b)$, then $\Upsilon ^nf(x^{[n]})\in {\bf K}(v,\phi )$ or ${\bar
{\Phi }}^nf(x^{(n)})\in {\bf K}(v,\phi )$, where $v$ is a marked
vector $v\in ({\bf K}^b)^{[n]}$ or $v\in ({\bf K}^b)^{(n)}$, where
$x\in {\bf K}^b$, $x^{[n]}\in ({\bf K}^b)^{[n]}$ or $x^{(n)}\in
({\bf K}^b)^{(n)}$ correspondingly.}
\par {\bf Proof.} Without loss of generality it can be assumed, that
the function $\phi $ is subadditive and increasing taking \par $\phi
_1(q) := \inf \{ \sum_{k=1}^n\phi (q_k): \quad \sum_{k=1}^nq_k \ge
q, q_k\ge 0 \} $, which is the largest increasing and subadditive
minorant of $\phi $. For the subadditive and increasing $\phi $
there is satisfied the inequality:
\par $(1)$ $\phi (q\epsilon )\le \phi ((1+[q])\epsilon )\le
(1+q)\phi (\epsilon )$ \\
for each $\epsilon >0$ and $q>0$, where $[q]$ denotes the integral
part of $q$ such that $[q]\le q$.
\par  If $S$ is a family of vectors such that it
spans ${\bf K}^b$ and $f$ belongs to ${\bf K}(u,\phi )$ for each
$u\in S$, then $f\in {\bf K}(\phi )$, since $b\in \bf N$ and
$|f(x+y)-f(x)|=|f(x+y)-f(x+y_2e_2+...+y_be_b)
+f(x+y_2e_2+...+y_be_b)-f(x+y_3e_3+...+y_be_b)+...+
f(x+y_be_b)-f(x)|\le \max (|f(x+y)-f(x+y_2e_2+...+y_be_b)|,
|f(x+y_2e_2+...+y_be_b)-f(x+y_3e_3+...+y_be_b)|,...,|f(x+y_be_b)-f(x)|)
\le C\max (\phi (|y_1|),...,\phi (|y_b|))\le C\phi (|y|)$ due to
increasing monotonicity of $\phi $ and the fact that $|y|=\max
(|y_1|,...,|y_b|)$, where $y=y_1e_1+...+y_be_b$, $y_1,...,y_b\in \bf
K$, $e_j=(0,...,0,1,0,...)\in {\bf K}^b$ with $1$ on the $j$-th
place and up to a $\bf K$-linear topological automorphism of ${\bf
K}^b$ onto itself we can choose such basis as belonging to $S$,
$j=1,...,b$, and $C>0$ is a constant.
\par Let us assume that for some point the statement of this lemma
is not true. We can suppose, that this is at $x^{[n]}=(0,...,0)\in
({\bf K}^b)^{[n]}$ or $x^{(n)}=0\in ({\bf K}^b)^{(n)}$ respectively
making a shift in a case of necessity. Then there exist sequences
$b_k>0$, $h_k\in \bf K$, $h_k\ne 0$, $\mbox{ }_kz^{[n]}\in ({\bf
K}^b)^{[n]}$ such that $\lim_{k\to \infty }b_k=\infty $, $\lim_{k\to
\infty }h_k=0$, $\lim_{k\to \infty }\mbox{ }_kz^{[n]}=0$ and
\par $(2)$ $|\Upsilon ^n f(\mbox{
}_kz^{[n]}+h_kv)-\Upsilon ^n f(\mbox{ }_kz^{[n]})|>b_k\phi (|h_k|)$
or
\par $(2')$ $|{\bar {\Phi }}^n f(\mbox{
}_kz^{(n)}+h_kv)-{\bar {\Phi }}^n f(\mbox{ }_kz^{(n)})|>b_k\phi
(|h_k|)$ \\
with $\lim_{k\to \infty } \mbox{ }_kz^{(n)}=0$ respectively, where
$0\ne v\in ({\bf K}^b)^{[n]}$ or $0\ne v\in ({\bf K}^b)^{(n)}$
correspondingly, $k=1,2,3,...$. Let the functions $u$ and $u_j$ be
as in the proof of Lemma 27. Choose $r_j\in \bf N$ such that
$|\mbox{ }_{r_j}z_0|\le |c_j|$, $\lim_{j\to \infty
}|c_j|^{n+1}b_{r_j}=\infty $, $|h_{r_j}|<|\pi c_jT_j|$. Thus $u\in
C^{\infty }({\bf K},{\bf K}^b)$. Now prove that at least for large
$j\in \bf N$ there is accomplished the inequality:
\par $(3)$  $|\Upsilon ^n(f\circ u)(\mbox{ }_jx^{[n]}+\mbox{ }_j\nu ^{[n]})-
\Upsilon ^n(f\circ u)(\mbox{ }_jx^{[n]})|>|\pi c_j^{n+1}| b_{r_j}
\phi (|\mbox{ }_j\nu ^{[n]}|) |\pi |^{l_0}$ or
\par $(3')$  $|{\bar {\Phi }}^n(f\circ u)(\mbox{ }_jx^{(n)}+
\mbox{ }_j\nu ^{(n)})- {\bar {\Phi }}^n(f\circ u)(\mbox{
}_jx^{(n)})|>|\pi c_j^{n+1}|
b_{r_j} \phi (|\mbox{ }_j\nu ^{[n]}|) |\pi |^{l_0}$, \\
where $|\mbox{ }_j\nu ^{[n]}|=|h_{r_j}v/c_j|$ or $|\mbox{ }_j\nu
^{(n)}|=|h_{r_j}v/c_j|$ with $c_j\ne 0$ for each $j$. Take without
loss of generality $|v|=1$. Together with the condition $\lim_{j\to
\infty }|c_j|^{n+1}b_{r_j}=\infty $ this will complete the proof. If
$|h|<|\pi T_j|$, then
\par $(4)$ $u_j(h)= \mbox{ }_{r_j}z_0+ c_j \sum_{k_1,k_2=0}^2
\sum_{i_1,i_2=1}^{b-1}\mbox{ }_{i_1,i_2}a_{k_1,k_2} h_{i_1}^{k_1}
h_{i_2}^{k_2}$ with $u_j(0)=\mbox{ }_{r_j}z_0$.\\
In formula $9(1)$ or $10(1)$ all terms with an amount of operators
$A_{j,v^{[k-1]},t_k}$ or $B_{j,v^{[k-1]},t_k}$ in it less than $n$
are in $C^{[1]}_{\phi }({\bf K},{\bf K})$. As in Lemma 27 reduce the
consideration to $\Upsilon ^nf(\mbox{ }_kz^{[n]})$ or ${\bar {\Phi
}}^nf(\mbox{ }_kz^{(n)})$ with prescribed fixed vectors $w_0,
w_1,...,w_{n-1}$ with $v^{[0]}=v^{[k],1}=w_0$ and $v^{[k],i}=w_l$
for each $2^{l-1}< i\le 2^l$, where $l=0,1,...,k$ and
$k=0,1,...,n-1$, vectors $v^{[k],i}\in {\bf K}^b$ are formed from
$v^{[k]}$ after excluding all zeros arising from $v_3^{[k]}=0$, or
$v=(w_0,...,w_0)$ with $\mbox{ }_kz^{(n)}= (\mbox{ }_kz_0;v;\mbox{
}_kt)$ respectively, where $\mbox{ }_kt=(\mbox{ }_kt_1,...,\mbox{
}_kt_n)\in {\bf K}^n$. In view of Lemma 16 we can consider the case
$(b-1,b)$ instead of $(1,b)$. Thus, from Formulas $(2,4)$ and
$12(2)$ it follows, that there exist expansion coefficients $\mbox{
}_{i_1,i_2}a_{k_1,k_2}\in {\bf K}^b$ with $|\mbox{
}_{i_1,i_2}a_{k_1,k_2}|\le 1$ for each $i_1, i_2, k_1, k_2$ and
there exists $j_0\in \bf N$, for which
\par $(5)$ $|\Upsilon ^n f\circ u(\mbox{ }_jx^{[n]}+\mbox{ }_j\nu ^{[n]})
- \Upsilon ^n f\circ u(\mbox{ }_jx^{[n]})| \ge |\pi
|^{l_0+s_0}|q_j|^n$\\ $ |\Upsilon ^n f( \mbox{ }_{r_j}z_0^{[n]}+
h_jv) - \Upsilon ^n f(\mbox{ }_{r_j}z_0^{[n]})| \ge
b_{r_j}|c_j^n|\phi (|c_j\mbox{ }_j\nu ^{[n]}|) |\pi |^{l_0}$ or
\par $(5')$ $|{\bar {\Phi }}^n
f\circ u(\mbox{ }_jx^{(n)}+\mbox{ }_j\nu ^{(n)})- {\bar {\Phi }}^n
f\circ u(\mbox{ }_jx^{(n)})| \ge |\pi |^{l_0+s_0}|q_j|^n$ \\ $|{\bar
{\Phi }}^n f( \mbox{ }_{r_j}z_0^{(n)}+ h_jv) - {\bar {\Phi }}^n
f(\mbox{ }_{r_j}z_0^{(n)})| \ge
b_{r_j}|c_j^n|\phi (|c_j\mbox{ }_j\nu ^{(n)}|) |\pi |^{l_0}$ \\
for each $j\ge j_0$, since $|\nu _j|<|\pi T_j|$, where $l_0\in \bf
N$ is a marked natural number, $\mbox{ }_j\tau _i$, $i=1,...,n$ are
parameters corresponding to $t_1,...,t_n$, but for the curve $u$
instead of $f$. There exists $j_0\in \bf N$ such that $|c_j|\le \min
(1,|\pi |^{-1}-1)$ for each $j>j_0$, where $\pi \in \bf K$, $0<|\pi
|<1$. In view of Formula $(1)$ for each $j>j_0$ we have $\phi
(|\mbox{ }_j\nu ^{[n]}|)\le (1+|c_j|^{-1})\phi (|c_j\mbox{ }_j\nu
^{[n]}|)\le |\pi c_j|^{-1} \phi (|c_j\mbox{ }_j\nu ^{[n]}|)$.
Therefore, the latter formula and Formula $(5)$ imply Formula $(3)$.
\par {\bf 31. Lemma.} {\it Let $f$ be a function $f: {\bf K}\to \bf
K$ such that $f(0)=0$ and $|f(t)|\le 1$ for each $t\in \bf K$ such
that $|t|\le |q|a$, where $q$ and $a$ are constants such that $q\in
\bf K$, $|q|>1$, $a>0$, and assume that
\par $(1)$ $|f(qt)-qf(t)|\le \max (b,C_1|t|^r)$ \\
for each $t\in \bf K$ with $|t|\le a$, where $0<r\le 1$, $b>0$ and
$C_1>0$ are constants. Then there exists a constant $C_2>0$ such
that
\par $(2)$ $|f(t)|\le \max (b,C_2|t|^r)$ \\
for each $t\in \bf K$ with $|t|\le |q|a$, where $C_2=\max
(a^{-r},|q|^{-1} a^rC_1|q|^{-r})$.}
\par {\bf Proof.} If $t\in \bf K$ is such that $a\le |t|\le |q|a$,
then Inequality $(2)$ is satisfied  with $C_2=a^{-r}$, since
$|f(t)|\le 1$ for such $t$. Now suppose that $0<|u|<a$ and
Inequality $(2)$ is satisfied for $t=qu$, then \par $|q| |f(t)|\le
\max (|f(qt)|,b,C_1|u|^r)\le \max (b,C_2|qu|^r,C_1|u|^r)$ \\  $=
\max (b, |u|^r\max (C_1,C_2|q|^r))$, hence \par $|f(t)|\le
|q|^{-1}\max (b,a^r\max (C_1,C_2|q|^r))\le \max (b,C_2|t|^r)$ \\ for
$C_2=\max (a^{-r},|q|^{-1} a^rC_1|q|^{-r})$, since $C_2|q|^r\ge
C_1$. On the other hand $B({\bf K},0,|q|a)\setminus \{ 0 \} $ is the
disjoint union of subsets $B({\bf K},0,|q|^ja)\setminus B({\bf
K},0,|q|^{j-1}a)$ for $j=1,0,-1,-2,...$. Therefore, proceeding by
induction by $j$ we get the statement of this lemma, since $f(0)=0$.
\par {\bf 32. Lemma.} {\it Let $\Omega $ be a finite set of vectors
$v\in {\bf K}^m$ which are pairwise $\bf K$-linearly independent,
$card (\Omega )\ge m$ and each subset of $\Omega $ consisting not
less than $m$ vectors has the $\bf K$-linear span coinciding with
${\bf K}^m$, where $m\ge 2$ is the integer. Suppose that for each
$v\in \Omega $ there is given a function $g_v: {\bf K}^m\to \bf K$
such that:
\par $(1)$ $|g_v(x)|\le 1$ for each $x\in {\bf K}^m$ with $|x|\le R$,
\par $(2)$ $|g_v(x+tv)-g_v(x)|\le |t|^r$ for each $x, x+tv\in {\bf K}^m$
with $|x|\le R$ and $|x+tv|\le R$, where $t\in \bf K$,
\par $(3)$ $|\sum_{v\in \Omega }(g_v(x)-g_v(y))|\le b$ for each
$|x|\le R$ and $|y|\le R$, where $R$ and $b$ are positive constants,
$0<r\le 1$. Then there exists a constant $C$, which may depend on
$(r,\Omega )$ such that
\par $(4)$ $|g_v(x)-g_v(y)|\le C\max (b,|x-y|^r)$ \\
for each $x, y\in {\bf K}^m$ such that $|x|\le R$, $|y|\le R$ and
each $v\in \Omega $.}
\par {\bf Proof.} Prove this lemma by induction on a number $n$ of
elements in $\Omega $. For $n=m=1$ Inequality $(4)$ is the
consequence of Inequality $(2)$. For $n=m\ge 2$ vectors
$v_1,...,v_m$ by the supposition of lemma are $\bf K$-linearly
independent. Then for each $x, y\in {\bf K}^m$ there exist $t_1,...,
t_m\in \bf K$ such that $x=y+t_1v_1+...+t_mv_m$. If $|x|\le R$ and
$|y|\le R$, then $B({\bf K}^m,0,R)=B({\bf K}^m,x,R)=B({\bf
K}^m,y,R)$ due to the ultrametric inequality. On the other hand,
$B({\bf K}^m,0,R)$ is the additive group, hence $y-x\in B({\bf
K}^m,0,R)$. Vectors $v_j$ have coordinates $v_j=(v_j^1,...,v_j^m)$,
where $v_j^k\in \bf K$, consequently, $|v_j|=\max
(|v_j^1|,...,|v_j^m|)$. Thus, $|x-y|= \max
(|t_1v_1^1+...+t_mv_m^1|,..., |t_1v_1^m+...+t_mv_m^m|)\le \max
(|t_1v_1|,...,|t_mv_m|)$. Choose $t_1,...,t_m$ such that
$|y+t_1v_1+...+t_kv_k|\le R$ also for each $k=1,...,m$. Therefore,
\par $|g_{v_j}(x)- g_{v_j}(y)|=|g_{v_j}(x)-g_{v_j}(y+t_1v_1+...
+t_{m-1}v_{m-1}) + g_{v_j}(y+t_1v_1+...+t_{m-1}v_{m-1})-...
-g_{v_j}(y+t_1v_1)+ g_{v_j}(y+t_1v_1) - g_{v_j}(y)|\le \max_{k=1}^m
|g_{v_j}(y+t_1v_1+...+t_kv_k)- g_{v_j}(y+t_1v_1+...+t_{k-1}v_{k-1})|
\le \max
(b,|t_1|^r,...,|t_m|^r)\le \max (b,|x-y|^r)$ \\
for each $|x|\le R$ and $|y|\le R$ as the consequence of
Inequalities $(2)$ and $(4)$ and the ultrametric inequality for each
$|x|\le R$ and $|y|\le R$, since $j=1,...,m$.
\par Further proceed by induction on $n$. From the preceding prove
it follows, that the statement of this lemma is true for $n=m$. Put
$\Omega = \Omega _0\cup \{ w \} $, where $w\notin \Omega _0$ and all
elements of $\Omega $ are pairwise linearly independent over the
field $\bf K$. Assume that the assertion of this lemma is true for
$\Omega _0$ and prove it for $\Omega $. For $v\in \Omega _0$ denote
$h_v(x,u) = h_v(x) = g_v(x+uw)-g_v(x)$, where $|x|\le R$, and
$|x+uv|\le R$. For these values of $x$ and $x+uv$ the function $h_v$
satisfies Conditions $(1,2)$. On the other hand, from $(2)$ for
$v=w$ and $(3)$ \par $|\sum_{v\in \Omega _0} (h_v(x)-h_v(y))|\le
\max (|\sum_{v\in \Omega }(g_v(x+uw)-g_v(x))|, |\sum_{v\in \Omega
}(g_v(y+uw)-g_v(y))|, |g_w((x+uw)-g_w(x)|, |g_w(y+uw)-g_w(y)|) \le
\max (b, |u|^r)$ \\
for each $|x|\le R$, $|y|\le R$ and $|uw|\le R$. Thus, $\{ h_v: v\in
\Omega _0 \} $ satisfies Condition $(3)$ with $\max (b,|u|^r)$
instead of $b$. By the induction hypothesis there exists
$C_1=const>0$, which may depend only on $\Omega _0$, $r$ and $R$
such that \par $(5)$ $|h_v(x)-h_v(y)|\le C_1\max (b,|u|^r,|x-y|^r)$ \\
for each $|x|\le R$, $|y|\le R$ and $|uw|\le R$ and $v\in \Omega
_0$. Take $y-x=(q-1)uw$ with $q\in \bf K$, $|q|>1$, hence $|q-1|>1$
and Inequality $(5)$ will take the form:
\par $(6)$ $|g_v(x+quw)-g_v(x+(q-1)uw) -g_v(x+uw) + g_v(x)|\le C_1\max
(b,|(q-1)u|^r,|(q-1)uw|^r)\le C_2\max (b,|u|^r)$, \\
when $|x|\le R$, $|quw|\le R$ and $v\in \Omega _0$, where $C_2\ge
C_1 |q-1|^r\max (1,|w|^r)$. Now set $s(u):=g_v(x+uw)-g_v(x)$ for
$v\in \Omega _0$ and from $(2)$ and $(6)$ and the ultrametric
inequality it follows, that
\par $|s(qu)-qs(u)|\le C_2\max (b,|u|^r)$, \\
when $|quw|\le R$. In view of Lemma 19 \par $(7)$ $|s(u)|=
|g_v(x+uw)-g_v(x)|\le C_3\max (b,|u|^r)$ \\
for each $|uw|\le R$, $|x|\le R$ and $v\in \Omega _0$, where
$C_3=\max (a^{-r},|q|^{-1} a^rC_2|q|^{-r})$. Interchanging roles of
$w$ and one of $v\in \Omega _0$ we obtain $(7)$ with $w$ in place of
$v$, that is, $(4)$ is proved for each $v\in \Omega $.
\par {\bf 33. Corollary.} {\it Let $v_1,...,v_n$ be pairwise $\bf
K$-linearly independent vectors in ${\bf K}^m$ and each subset
consisting not less than $m$ of these vectors has the $\bf K$-linear
span coinciding with ${\bf K}^m$ and let $g_k$ be locally bounded
functions from ${\bf K}^m$ into $\bf K$, $0<r\le 1$. If $g_k\in
Lip(v_k,r)$ for each $k$ and $\sum_{k=1}^ng_k(x)=0$ identically by
$x\in {\bf K}^m$, then $g_k\in Lip (r)$ for each $k$.}
\par {\bf Proof.} If $c\in \bf K$, $c\ne 0$ is small enough, then
the functions $cg_k$ satisfy assumptions of Lemma 32 with $b=0$.
\par {\bf 34. Remark.} We can mention, that apart from the classical
case over $\bf R$ this lemma is true also for $r=1$ due to the
ultrametric inequality, which is stronger than the usual triangle
inequality.
\par {\bf 35. Definition.} Let $v\in {\bf K}^b$ and $v\ne 0$.
We say that a function $f: {\bf K}^b\to \bf K$ is continuous in the
direction $v$ if $f(x+tv)$ converges to $f(x)$ uniformly by $x$ on
bounded closed sets as $t$ tends to zero.
\par Mention that in a particular case of a locally compact field
$\bf K$ a bounded closed subset is compact.
\par {\bf 36. Lemma.} {\it Suppose that $f\in C^0({\bf K}^b,
{\bf K})$ and $\Upsilon ^1f(x,w,t)$ is continuous or uniformly
continuous on $V^{[1]}$ in the direction $v^{[1]}$ with
$v_2^{[1]}\ne 0$ and $v_3^{[1]}\ne 0$, where $V^{[1]} := \{
(x,v,t)\in U^{[1]}: |v|=1 \} $, $U$ is open in ${\bf K}^b$. Then
$\Upsilon ^1f(x,v_2^{[1]},t)$ is continuous or uniformly continuous
by $(x,t)$, $(x,v,t)\in U^{[1]}$ or $(x,v,t)\in V^{[1]}$
respectively.}
\par {\bf Proof.} Assume the contrary, that
$\Upsilon ^1f(x,v_2^{[1]},t)$ is not continuous by $(x,t)$. Making a
shift in case of necessity we can suppose that $\Upsilon
^1f(x,v_2^{[1]},t)$ is not continuous by $(x,t)$ at $0$ or is not
uniformly continuous on $V^{[1]}$. Therefore, there exists a
sequence $ \{ x_n^{[1]}\in ({\bf K}^b)^{[1]}: n\in {\bf N} \} $ such
that $|\Upsilon ^1f(x_n^{[1]}) - \Upsilon ^1f(0)|>\epsilon $ for
each $n$ or with $x_0^{[1]}\in V^{[1]}$ instead of $0$ and a family
of sequences parametrized by $x_0^{[1]}$ and $\sup_{x_0^{[1]}\in
V^{[1]}}|\Upsilon ^1f(x_n^{[1]}) - \Upsilon ^1f(x_0^{[1]})|>\epsilon
$ correspondingly, where $\epsilon
>0$ is a constant, $x_n^{[1]}=(x_n,v_2^{[1]},t_n)$, $\lim_{n\to
\infty }(x_n,t_n)=x_0^{[1]}$. But in accordance with Definition 21
there exists $\delta
>0$ independent of $n$ such that $|\Upsilon ^1f(x_n^{[1]}+v^{[1]}\tau )-
\Upsilon ^1f(v^{[1]}\tau )|>\epsilon |\pi |$ or $\sup_{x_0^{[1]}\in
V^{[1]}} |\Upsilon ^1f(x_n^{[1]}+v^{[1]}\tau )- \Upsilon
^1f(v^{[1]}\tau )|>\epsilon |\pi |$ for each $n$ and each $\tau \in
\bf K$ with $|\tau |\le \delta $. On the other hand, $\Upsilon
^1f(x_n+v_1^{[1]}\tau ,w_n+v_2^{[1]}\tau ,t_n+v_3^{[1]}\tau )-
\Upsilon ^1f(v^{[1]}\tau )= [f(x_n+v_1^{[1]}\tau +(w_n+v_2^{[1]}\tau
)(t_n+v_3^{[1]}\tau ))-f(x_n+v_1^{[1]}\tau )]/(t_n+v_3^{[1]}\tau ) -
[f(v_1^{[1]}\tau +v_2^{[1]}\tau v_3^{[1]}\tau )- f(v_1^{[1]}\tau
)]/(v_3^{[1]}\tau )$, where $w_n=v_2^{[1]}$. But
\par $\lim_{n\to \infty }[f(x_n+v_1^{[1]}\tau +(w_n+v_2^{[1]}\tau )
(t_n+v_3^{[1]}\tau ))/(t_n+v_3^{[1]}\tau ) - f(v_1^{[1]}\tau
)/(v_3^{[1]}\tau )]=0$ and
\par $\lim_{n\to \infty } [f(x_n+v_1^{[1]}\tau
)/(t_n+v_3^{[1]}\tau ) - f(v_1^{[1]}\tau )]/(v_3^{[1]}\tau )=0$ \\
for $v_3^{[1]}\tau \ne 0$ pointwise or uniformly respectively. If
$v_3^{[1]}\tau =0$, then $\Upsilon ^1f(x_n+v_1^{[1]}\tau
,w_n+v_2^{[1]}\tau ,t_n+v_3^{[1]}\tau )- \Upsilon ^1f(v^{[1]}\tau
)=\Upsilon ^1f(x_n+v_1^{[1]}\tau ,w_n+v_2^{[1]}\tau ,0)-\Upsilon
^1f(v_1^{[1]}\tau ,w_n+v_2^{[1]}\tau ,0)$, but the latter difference
tends to zero as $\tau $ tends to zero uniformly by $n$ or also
uniformly by the family of sequences parametrized by $x_0^{[1]}$
respectively in accordance with the supposition of lemma. Thus we
get the contradiction with our supposition, hence $\Upsilon
^1f(x,v_2^{[1]},t)$ is continuous or uniformly continuous by $(x,t)$
correspondingly.
\par {\bf 37. Definition.} Denote by either $C^{[n]}_{\phi ,b}(U,Y)$ or
$C^n_{\phi ,b}(U,Y)$ spaces of all functions $f\in C^{[n]}_{\phi
}(U,Y)$ or $f\in C^n_{\phi }(U,Y)$ such that $f^{[k]}(x^{[k]})$ or
${\bar {\Phi }}^kf(x^{(k)})$ is uniformly continuous on a subset
either $V^{[k]} := \{ x^{[k]}\in U^{[k]}: |v^{[q]}_1|=1; |\mbox{
}_lv^{[q]}_2t_{q+1}|\le 1, |v^{[q]}_3|\le 1 \quad \forall l, q \} $
or $V^{(k)} := \{ x^{(k)}\in U^{(k)}: |v_j|=1 \quad \forall j \} $
for each $k=0,1,...,n$ with finite norms either \par $\|
f\|_{[n]}:=\| f\|_{[n],\phi }:=\max (C, \sup_{k=0,...,n; x^{[k]}\in
V^{[k]}} |f^{[k]}(x^{[k]})|)$ or
\par $\| f\|_{n}:=\| f\|_{n,\phi }:=
\max (C, \sup_{k=0,...,n; x^{(k)}\in V^{(k)}}
|{\bar {\Phi }}^kf(x^{(k)})|)$, \\
where $0\le C<\infty $ is the least constant satisfying $29(1)$ or
$29(2)$ for $\Upsilon ^nf(x^{[n]})$ or ${\bar {\Phi }}^nf(x^{(n)})$
respectively instead of $f$. For $\phi (q)=q^r$ we denote
$C^{[n]}_{\phi }(U,Y)$ by $C^{[n],r}(U,Y)$ and $C^n_{\phi ,b}(U,Y)$
by $C^{n,r}_b(U,Y)$, $0\le r\le 1$.
 For $r=0$ we put $C^{[n],0}=C^{[n]}$,
$C^{[n],0}_b=C^{[n]}_b$ and $C^{n,0}=C^{n}$, $C^{n,0}_b=C^{n}_b$,
with $C=0$ in the definition of the norm.
 As usually $C^{[\infty
]}(U,Y):=\bigcap_{k=0}^{\infty }C^{[k]}(U,Y)$ and $C^{\infty
}(U,Y):=\bigcap_{k=0}^{\infty }C^k(U,Y)$ and $C^{[\infty
]}_b(U,Y):=\bigcap_{k=0}^{\infty }C^{[k]}_b(U,Y)$ and $C^{\infty
}(U,Y):=\bigcap_{k=0}^{\infty }C^k_b(U,Y)$, where the topology of
the latter two spaces is given by the family of the corresponding
norms.
\par In the case of a locally compact field $\bf K$ and a compact
clopen (closed and open at the same time) domain $U$ we have
$C^{[k]}_b(U,Y)=C^{[k]}(U,Y)$ and $C^{k}_b(U,Y)=C^{k}(U,Y)$, though
for non locally compact $\bf K$ they are different $\bf K$-linear
spaces.
\par {\bf 38. Theorem.} {\it Suppose that $f: {\bf K}^m\to
\bf K$, $m\in \bf N$ and $f\circ u\in C^s_{\phi }({\bf K},{\bf K})$
or $f\circ u\in C^s_{\phi ,b}({\bf K},{\bf K})$ or $f\circ u\in
C^{[s]}_{\phi }({\bf K},{\bf K})$ or $f\circ u\in C^{[s]}_{\phi
,b}({\bf K},{\bf K})$ for each $u\in C^{\infty }({\bf K},{\bf K}^m)$
or $u\in C^{\infty }_b({\bf K},{\bf K}^m)$ or $u\in C^{[\infty ]
}({\bf K},{\bf K}^m)$ or $u\in C^{[\infty ]}_b({\bf K},{\bf K}^m)$,
where $s$ is a nonnegative integer, $\phi : (0,\infty )\to (0,\infty
)$ such that $\lim_{y\to 0}\phi (y)=0$, then $f\in C^s({\bf
K}^m,{\bf K})$ or $f\in C^s_b({\bf K}^m,{\bf K})$ or $f\in
C^{[s]}({\bf K}^m,{\bf K})$ or $f\in C^{[s]}_b({\bf K}^m,{\bf K})$
respectively.}
\par {\bf Proof.} In view of Lemma 21 it is sufficient to
prove that \\ ${\bar {\Phi
}}^nf(x;e_{j(1)},...,e_{j(n)};t_1,...,t_n)$ is in
$C^0(U^{(n)}_{j(1),...,j(n)},Y)$ or
$C^0_b(V^{(n)}_{j(1),...,j(n)},Y)$ or $\Upsilon ^nf(x^{[n]})$ is in
$C^0(U^{[n]}_{j(0),...,j(n)},Y)$ or
$C^0_b(V^{[n]}_{j(0),...,j(n)},Y)$ respectively for each
$n=1,2,...,s$ and each $j(1),...,j(n)\in \{ 1,...,m \} $ or $j(i)\in
\{ 1,...,m(i) \} $, $i=0, 1,..., n$. If $\Upsilon ^{m+1}f$ or ${\bar
{\Phi }}^{m+1}f$ is locally bounded, then $\Upsilon ^mf$ or ${\bar
{\Phi }}^mf$ is continuous respectively. Applying Lemma 27 by
induction we get that $\Upsilon ^nf(x^{[n]})$ is a locally bounded
function on $({\bf K}^m)^{[n]}$ and ${\bar {\Phi }}^nf(x^{(n)})$ is
a locally bounded function on $({\bf K}^m)^{(n)}$. In view of Lemma
30 each $\Upsilon ^kf(x^{[k]})$ or ${\bar {\Phi }}^kf(x^{(k)})$ is
continuous in each direction $v$ for each $k=1,...,s$, where $v\in
({\bf K}^m)^{[k]}$ or $v\in ({\bf K}^m)^{(k)}$ correspondingly. On
the other hand, by induction on $k$ we have that in accordance with
Lemma 36 $\Upsilon ^kf(x^{[k]})$ or ${\bar {\Phi
}}^kf(x;e_{j(1)},...,e_{j(k)};t_1,...,t_k)$ is continuous on
$U^{[k]}_{j(0),...,j(k)}$ or $U^{(k)}_{j(1),...,j(k)}$ or bounded
uniformly continuous respectively on $V^{[k]}_{j(0),...,j(k)}$ or
$V^{(k)}_{j(1),...,j(k)}$  for bounded $U$ for each
$j(1),...,j(k)\in \{ 1,...,m \}$ or $j(i)\in \{ 1,...,m(i) \} $ for
all $i=0,1,...,k$.
\par {\bf 39. Theorem.} {\it Let $f: {\bf K}^m\to {\bf K}^n$,
$m, n\in \bf N$. Let also $f\circ u\in C^{s,r}({\bf K},{\bf K}^n)$
or $f\circ u\in C^{s,r}_b({\bf K},{\bf K}^n)$ or $C^{[s],r}({\bf
K},{\bf K}^n)$ or $C^{[s],r}_b({\bf K},{\bf K}^n)$ for each $u\in
C^{\infty }({\bf K},{\bf K}^m)$ or $u\in C^{\infty }_b({\bf K},{\bf
K}^m)$ or $C^{[\infty ]}({\bf K},{\bf K}^m)$ or $C^{[\infty
]}_b({\bf K},{\bf K}^m)$ correspondingly, where $s$ is a nonnegative
integer, $0\le r\le 1$, then $f\in C^{s,r}({\bf K}^m,{\bf K}^n)$ or
$f\in C^{s,r}_b({\bf K}^m,{\bf K}^n)$ or $C^{[s],r}({\bf K}^m,{\bf
K}^n)$ or $C^{[s],r}_b({\bf K}^m,{\bf K}^n)$ respectively.}
\par {\bf Proof.} If $s=0$ and $0\le r\le 1$, then the assertion of this
theorem follows from Lemmas 27 and 30. For $r>0$ by Theorem 38 $f\in
C^s({\bf K}^m,{\bf K}^n)$ or $f\in C^s_b({\bf K}^m,{\bf K}^n)$ or
$C^{[s]}({\bf K}^m,{\bf K}^n)$ or $C^{[s]}_b({\bf K}^m,{\bf K}^n)$
respectively. From Lemma 21 we infer, that it is sufficient to prove
that ${\bar {\Phi }}^nf(x;e_{j(1)},...,e_{j(n)};t_1,...,t_n)$ is in
$C^{0,r}(U^{(n)}_{j(1),...,j(n)},Y)$ or
$C^{0,r}_b(V^{(n)}_{j(1),...,j(n)},Y)$ or $\Upsilon ^nf(x^{[n]})\in
C^{0,r}(U^{[n]}_{j(0),...,j(n)},Y)$ or
$C^{0,r}_b(V^{[n]}_{j(0),...,j(n)},Y)$ respectively for each
$n=1,2,...,s$ and each $j(1),...,j(n)\in \{ 1,...,m \} $, $j(i)\in
\{ 1,...,m(i) \} $, $i=0,1,...,n$. Prove this by induction by $n$.
For $n=0$ it was proved above. Let it be true for $n=0,...,k$ and
prove it for $n=k+1\le s$. For this consider Formula $10(1)$ or
$9(1)$. On the right hand side of it all terms having a total degree
of $f$ by operators $B$ or $A$ less than $k+1$ are in
$C^{0,r}(U^{(n)},Y)$ or $C^{0,r}_b(V^{(n)},Y)$ or
$C^{0,r}(U^{[n]},Y)$ or $C^{0,r}_b(V^{(n)},Y)$ respectively by the
induction hypothesis, since $u\in C^{\infty }({\bf K},{\bf K}^m)$ or
$u\in C^{\infty }_b({\bf K},{\bf K}^m)$ or $C^{[\infty ]}({\bf
K},{\bf K}^m)$ or $C^{[\infty ]}_b({\bf K},{\bf K}^m)$
correspondingly. Therefore, it remains to prove, that the sum
\par $(i)$ $[\sum_{j_1,...,j_n}(B_{j_n,v^{(n-1)},t_n}...
B_{j_1,v^{(0)},t_1} f\circ u) ({\bar {\Phi }}^1\circ
p_{j_n}{\hat S}_{j_{n-1}+1,v^{(n-2)}t_{n-1}}$ \\
$...{\hat S}_{j_1+1,v^{(0)}t_1}u^{n-1}) (P_n{\bar {\Phi }}^1\circ
p_{j_{n-1}}{\hat S}_{j_{n-2}+1,v^{(n-3)},t_{n-2}}... {\hat
S}_{j_1+1,v^{(0)}t_1}u^{n-2})...(P_n...P_2{\bar {\Phi }}^1\circ
p_{j_1}u)]$ \\ is in $C^{0,r}(U^{(n)},Y)$ or $C^{0,r}_b(V^{(n)},Y)$
or corresponding sum by compositions of $A_{j_k,v^{[k-1]},t_k}$ is
in $C^{0,r}(U^{[n]},Y)$ or $C^{0,r}_b(V^{[n]},Y)$ respectively. In
accordance with the proof above it is sufficient to demonstrate this
for $v^{(n-1)}=(e_{j(1)},...,e_{j(n)})$ for each $j(1),...,j(n)\in
\{ 1,...,m \} $ or $v^{[i]}=e_{j(i)}$ with $j(i)\in \{ 1,...,m(i) \}
$ and $i=0,1,...,n$, where $v_0^{(l)}=v_{l+1}=e_{j(l+1)}$,
$l=0,...,n-1$. By the induction hypothesis ${\bar {\Phi
}}^lf(x;e_{j(1)},...,e_{j(l)};t_1,...,t_l)$ is in
$C^{0,r}(U^{(l)}_{j(1),...,j(l)},Y)$ or
$C^{0,r}_b(V^{(l)}_{j(1),...,j(l)},Y)$ or $\Upsilon ^lf(x^{[l]})\in
C^{0,r}(U^{[l]}_{j(0),...,j(l)},Y)$ or
$C^{0,r}_b(V^{[l]}_{j(0),...,j(l)},Y)$ respectively for each
$l=1,2,...,k$ and each $j(1),...,j(l)\in \{ 1,...,m \} $, $j(i)\in
\{ 1,...,m(i) \} $, $i=0,...,n$. In view of Corollary 18 and Lemma
32 functions ${\bar {\Phi
}}^nf(x;e_{j(1)},...,e_{j(n)};t_1,...,t_n)$ or $\Upsilon
^nf(x^{[n]})$ belong to $Lip (v,r)$ by $(x,t_1,...,t_n)$ or
$x^{[n]}\in U^{[n]}_{j(0),...,j(n)}$, where $v=(e_{j(n)};l_k)\in
{\bf K}^{m+n}$, $e_j\in {\bf K}^m$, $l_k=(0,...,0,1,0,...,0)\in {\bf
K}^n$ with $1$ on the $k$-th place, or
$v=(e_{j(0)},...,e_{j(n)};l_k)$ with $e_{j(i)}\in {\bf K}^{m(i)}$
respectively. By Corollary 18 each ${\bar {\Phi
}}^nf(x;e_{j(1)},...,e_{j(n)};t_1,...,t_n)$ or $\Upsilon
^nf(x^{[n]})|_{U^{[n]}_{j(0),...,j(n)}}$ belongs to $Lip (r)$ by
$(x,t_1,...,t_n)$ or in addition bounded uniformly lipschitzian on
$V^{(n)}_{j(1),...,j(n)}$ or $V^{[n]}_{j(0),...,j(n)}$ respectively.
In accordance with Lemma 21 this proves the theorem.
\par {\bf 40. Theorem.} {\it Let $f: {\bf K}^m\to {\bf K}^l$,
$m\in \bf N$. Suppose also that $f\circ u\in C^{\infty }({\bf
K},{\bf K})$ or $f\circ u\in C^{\infty }_b({\bf K},{\bf K})$ or
$C^{[\infty ]}({\bf K},{\bf K})$ or $C^{[\infty ]}_b({\bf K},{\bf
K})$ for each $u\in C^{\infty }({\bf K},{\bf K}^m)$ or $u\in
C^{\infty }_b({\bf K},{\bf K}^m)$ or $C^{[\infty ]}({\bf K},{\bf
K}^m)$ or $C^{[\infty ]}_b({\bf K},{\bf K}^m)$, then $f\in C^{\infty
}({\bf K}^m,{\bf K}^l)$ or $f\in C^{\infty }_b({\bf K}^m,{\bf K}^l)$
or $C^{[\infty ]}({\bf K}^m,{\bf K}^l)$ or $C^{[\infty ]}_b({\bf
K}^m,{\bf K}^l)$ respectively.}
\par {\bf Proof.} Apply either Theorem 39 for each $s\in \bf N$ and
$r=0$ or Theorem 38 for each $s\in \bf N$ and $\phi (q)=q^r$ with
$0<r<1$, since $C^{s,r}(U,Y)\subset C^{s+1}(U,Y)$ and
$C^{s,r}_b(U,Y)\subset C^{s+1}_b(U,Y)$ and $C^{\infty
}(U,Y):=\bigcap_{n=0}^{\infty }C^n(U,Y)$ and $C^{\infty
}_b(U,Y):=\bigcap_{n=0}^{\infty }C^n_b(U,Y)$ and $C^{[\infty
]}(U,Y):=\bigcap_{n=0}^{\infty }C^{[n]}(U,Y)$ and $C^{[\infty
]}_b(U,Y):=\bigcap_{n=0}^{\infty }C^{[n]}_b(U,Y)$.
\par {\bf 41. Theorem.} {\it Let $h_j(y)$ be $C^{\infty }({\bf K},
{\bf K})$ functions such that \par $(1)$ $h_j(0)=0$ for each
$j=0,1,...,m$,
\par $(2)$ $\lim_{0\ne y\to 0} h_{j-1}(y)/h_j(y)^n=0$ for each $n\in \bf
N$ and $j=1,...,m$,
\par $(3)$ $\lim_{0\ne y\to 0}h_m(y)/y^n=0$ for every $n\in \bf N$.
Put $h(y)=(h_1(y),...,h_m(y))$ and suppose that $g\in C^{\infty }
({\bf K}^b,{\bf K})$ is not identically zero and $g(x)=0$ for each
$|x|>1$. Define $f: {\bf K}^{m+1}\to \bf K$ by the formula: \par
$(4)$ $f(x,y)=g((x-h(y))/h_0(y))$ for each $y\ne 0$ and $x\in {\bf
K}^m$, $f(x,0)=0$ for each $x$. Then $f\circ u\in C^{\infty }({\bf
K}^m,{\bf K})$ for each locally analytic function $u: {\bf K}^m\to
{\bf K}^{m+1}$, but $f$ is discontinuous.}
\par {\bf Proof.} First demonstrate that $f$ is not continuous at
$(0,0)$. We have $f(x,0)=0$. But take a sequence $(x_n,y_n)$ such
that $\lim_{n\to \infty }(x_n,y_n)=0$ with  $\epsilon \le
|(x_n-h(y_n))/h_0(y_n)|\le 1$ and $|g(z_n)|\ge \delta $ with
$z_n:=(x_n-h(y_n))/h_0(y_n)$, which is possible since $\lim_{0\ne
y\to 0} |h(y)/h_0(y)|=\infty $ and $g$ is continuous and non zero,
where $\epsilon >0$ and $\delta >0$ are constants. For this sequence
we have $|f(x_n,y_n)|\ge \delta $ for each $n$. But for the sequence
$(x_n,y_n)$ such that $|(x_n-h(y_n))/h_0(y_n)|> 1$ we have
$f(x_n,y_n)=0$, since $g(z_n)=0$ for $|z_n|>1$. Thus $f$ is
discontinuous at $(0,0)$.
\par Now take a locally analytic function $u: {\bf K}^m\to
{\bf K}^{m+1}$ and consider the composition $f\circ u$. Take a
nontrivial analytic function $w(x,y)$ from a neighborhood of zero in
${\bf K}^{m+1}$ into $\bf K$ such that $w\circ u(y)=0$ in a
neighborhood of $y_0\in {\bf K}^m$, where $u(y_0)=0$. Prove that for
functions $h_j(y)$ satisfying Conditions $(1-3)$ there exist
constants $C>0$ and $\delta >0$ such that
\par $(5)$ $|w(x,y)|\ge C |h_0(y)|^n$, \\
when $|x-h(y)|\le |h_0(y)|$, $0<|y|<\delta $. If prove $(5)$, then
from $(x,y)=u(t)$ for $t$ in a neighborhood of $y_0$ it follows,
that $w(x,y)=0$ and by $(5)$ we have that $|x-h(y)|>|h_0(y)|$, hence
$f(x,y)=0$.
\par Consider an analytic function $q$ from a neighborhood of zero
in ${\bf K}^l$ into $\bf K$. Then we can write it in the form: \par
$(6)$ $q(x)=x_1^ks(x_2,...,x_l)+x_1^{k+1}r(x)$, \\
where $s$ and $r$ are analytic functions and $s$ is not identically
zero, $1\le l\le m$, $x_1,...,x_l\in \bf K$. Then $|q(h(y))|\ge
|h_1(y)|^k(|s(h_2(y),...,h_l(y))|-C|h_1(y)|)\ge
|h_1(y)|^k(C|h_2(y)|^n-C|h_1(y)|)\ge C|\pi h_1(y)|^{k+n}$ \\
for each $y$ with $0<|y|<\delta $ with suitable $\delta >0$ and
$n\in \bf N$, where $\pi \in \bf K$, $0<|\pi |<1$. Then induction by
$l$ gives from $(6)$ that for an arbitrary nontrivial analytic
function $q$ from a neighborhood of zero in ${\bf K}^m$ into $\bf K$
there exist constants $C>0$ and $\delta
>0$ and $n\in \bf N$ such that
\par $(7)$ $|q(h(y))|\ge C|h_1(y)|^n$ for each $0<|y|<\delta $. \\
In particular, from $(7)$ it follows, that there exist $C>0$,
$\delta >0$ and $n\in \bf N$ such that
\par $(8)$ $|w(h(y),y)|\ge C|h_1(y)|^n$ for each $0<|y|<\delta $.
\par It remains to show that $(8)$ implies $(5)$.
Take $C>0$ so large that $|grad_x w(x,y)|\le C$ in some neighborhood
of zero and assume that $|x-h(y)|\le |h_0(y)|$. Then there exists
$\delta >0$ such that $|w(x,y)|\ge |w(h(y),y)|-C|x-h(y)|\ge
C|h_1(y)|^n-C|h_0(y)|\ge C|\pi h_1(y)|^n$ for each $0<|y|<\delta $.
Thus $f\circ u=0$ in a neighborhood of each point $y_0\in {\bf K}^m$
such that $u(y_0)=0$.
\par For example, we can take either $h_j(y)=\sum_n a_n
\pi ^{n^2(m-j+1)+n}$ for each $y=\sum_n a_n \pi ^n\in \bf K$, where
$a_n\in \bf K$ belong to the finite set of representatives of
distinct classes in the finite factor field $B({\bf K},0,1)/B({\bf
K},0,|\pi |)$, $\pi \in \bf K$, $0<|\pi |<1$, $\bf K$ is a locally
compact field of zero characteristic and $|\pi |$ is the largest
generator among those less than one of the valuation group $\Gamma
_{\bf K}$ of $\bf K$, or $h_j(y)=\sum_n a_n \theta ^{n^2(m-j+1)+n}$
for each $y=\sum_n a_n\theta ^n\in {\bf F}_{p^k}(\theta )$, where
$a_n\in {\bf F}_{p^k}$, $p$ is a prime number, $k\in \bf N$, ${\bf
F}_{p^k}(\theta )$ is a locally compact field of characteristic
$char ({\bf F}_{p^k}(\theta ))=p>0$, ${\bf F}_{p^k}$ is a finite
field of $p^k$ elements.
\par {\bf 42. Theorem.} {\it There exists a discontinuous function
$f: {\bf K}^m\to \bf K$ such that $f\circ u\in C^{\infty }({\bf
K}^{m-1},{\bf K})$  for each locally analytic function $u: {\bf
K}^{m-1}\to {\bf K}^b$, where $m\ge 2$.}
\par {\bf Proof.} This theorem follows from Theorem 41. Another
its proof is the following. Let $f\in C^{\infty }({\bf K}^2\setminus
\{ 0 \}, {\bf K})$ and let $f$ be non constant with $f(x_1,x_2)=0$,
when $x_1x_2=0$. For simplicity let ${\bf K}$ be a locally compact
field of zero characteristic. Take an analytic function $g: {\bf K}
\to \bf K$ such that $\lim_{|y|\to \infty } g(y)=0$. Such functions
exist due to Example 43.1 of Section 43 in \cite{sch1}. Moreover,
they can be chosen such that $|g(y)|\le \epsilon _j$ for $|y|=|\pi
|^{-j}$ for each $j=0,1,2,...$ and a sequence $\{ \epsilon _j>0: j
\} $, which in particular may also tend to zero. Then consider the
function $g(1/x_2)$ and put $h(x_1,x_2):=f(x_1,g(1/x_2))$, where $f$
is homogeneous of degree zero. Since $f\in C^{\infty }({\bf
K}^2\setminus \{ 0 \}, {\bf K})$ it remains to show that $f\circ
u\in C^{\infty }$ in a neighborhood of $y=0$, if $u(0)=(0,0)$. If
$u_1$ coincides with zero, then $h$ is identically zero. If
$u_1(0)=0$ and $u_1$ is not identically zero, then due to
analyticity there exists $k\in \{ 1,2,... \} $ such that
$u_1(t)=t^kv_1(t)$ and $v_1$ is locally analytic and $v_1(0)\ne 0$.
From $u_2(0)=0$ it follows that $g(1/u_2(t))=t^kv_2(t)$, where $v_2$
is locally analytic and $v_2(0)=0$. We can take, for example,
$\epsilon _j= |\pi |^{j^2}$. Since $f(x_1,0)=0$ and $f$ is
homogeneous of degree zero, then
$h(u_1(t),u_2(t))=f(t^kv_1(t),t^kv_2(t))=f(v_1(t),v_2(t))$ for each
$t\in \bf K$. Since $v_1(0)\ne 0$, then $f\circ u\in C^{\infty }$ in
a neighborhood of zero.
\par {\bf 43. Remark.} In the non archimedean case analogs of
classical theorems over $\bf R$ such as 3 and 10 \cite{boman} are
not true due to the ultrametric inequality for the non archimedean
norm, and since if a function $f$ is homogeneous, then ${\bar {\Phi
}}^k$ need not be homogeneous for $k\ge 1$. Theorem 2 from
\cite{boman} in the non archimedean case is true in the stronger
form due to the ultrametric inequality (see Theorem 38 above). The
notion of quasi analyticity used in the classical case in
\cite{boman} has not sense in the non archimedean case because of
the necessity to operate with ${\bar {\Phi }}^kf$ instead of $D^kf$.
It leads naturally to the local analyticity in the non archimedean
case. In the latter case the exponential function has finite radius
of convergence on $\bf K$ with $char ({\bf K})=0$. Therefore, in the
proof of Theorem 40 it was used specific feature of the non
archimedean analysis of analytic functions for which an analog of
the Louiville theorem is not true (see also \cite{sch1}).
\par Using the particular variant of Theorem 38 with $s=r=0$ it is
easy to prove the following theorem.
\par {\bf Theorem.} {\it Let $f: {\bf K}^m\to {\bf K}^l$, $f\circ
u\in C^n({\bf K}^2,{\bf K}^l)$ or $f\circ u\in C^n_b({\bf K}^2,{\bf
K}^l)$ or $C^{[n]}({\bf K}^2,{\bf K}^l)$ or $C^{[n]}_b({\bf
K}^2,{\bf K}^l)$ for each $u\in C^{\infty }({\bf K}^2,{\bf K}^m)$ or
$u\in C^{\infty }_b({\bf K}^2,{\bf K}^m)$ or $C^{[\infty ]}({\bf
K}^2,{\bf K}^m)$ or $C^{[\infty ]}_b({\bf K}^2,{\bf K}^m)$, where
$m\ge 2$ and $n\ge 1$. Then $f\in C^n({\bf K}^m,{\bf K}^l)$ or $f\in
C^n_b({\bf K}^m,{\bf K}^l)$ or $C^{[n]}({\bf K}^m,{\bf K}^l)$ or
$C^{[n]}_b({\bf K}^m,{\bf K}^l)$ correspondingly.}
\par {\bf Proof.} Put $u(y)=\sum_{j=1}^my_1^je_j+w(y_2)$,
where $y=(y_1,y_2)\in {\bf K}^2$, $e_j\in {\bf K}^m$, $w\in
C^{\infty }({\bf K},{\bf K}^m)$ or $C^{\infty }_b({\bf K},{\bf
K}^m)$ or $C^{[\infty ]}({\bf K},{\bf K}^m)$ or $C^{[\infty
]}_b({\bf K},{\bf K}^m)$. Therefore, $u\in C^{\infty }({\bf
K}^2,{\bf K}^m)$ or $C^{\infty }_b({\bf K}^2,{\bf K}^m)$ or
$C^{[\infty ]}({\bf K}^2,{\bf K}^m)$ or $C^{[\infty ]}_b({\bf
K}^2,{\bf K}^m)$. In view of Formula $10(1)$ or $9(1)$ and Lemmas
11, 12 or Corollary 14 for ${\bar {\Phi }}^kf\circ u(y^{(k)})$ or
$\Upsilon ^nf\circ u$ by induction we get that each ${\bar {\Phi
}}^kf(w(y_2),e_{j(1)},...,e_{j(k)};t_1,...,t_k)$ or $\Upsilon
^nf(x^{[n]})|_{V^{[n]}_{j(0),...,j(n)}}$ with $x=w(y_2)$ is
continuous or uniformly continuous. Therefore, from Theorem 39 with
$s=r=0$ it follows, that each ${\bar {\Phi
}}^kf(x;e_{j(1)},...,e_{j(k)};t_1,...,t_k)$ or $\Upsilon
^nf(x^{[k]})$ is continuous on $U^{(k)}_{j(1),...,j(k)}$ or
$U^{[n]}_{j(0),...,j(n)}$ or uniformly continuous on
$V^{(k)}_{j(1),...,j(k)}$ or $V^{[k]}_{j(0),...,j(k)}$ respectively
for each $k=1,...,n$, hence by Lemma 21 $f\in C^n({\bf K}^m,{\bf
K}^l)$ or $f\in C^n_b({\bf K}^m,{\bf K}^l)$ or $C^{[n]}({\bf
K}^m,{\bf K}^l)$ or $C^{[n]}_b({\bf K}^m,{\bf K}^l)$
correspondingly.

\par {\underline {Acknowledgement.}} The author is sincerely
grateful to Prof. K.-H. Neeb and Dr. H. Gl\"ockner for hospitality
and discussions of this matter at Mathematical Institute of
Darmstadt University of Technology and DFG for support.

\end{document}